\documentclass[11pt]{article}

\usepackage[margin=1in]{geometry}
\usepackage{setspace}
\usepackage{amsmath,amssymb,amsthm}
\usepackage{bm}

\usepackage{authblk}

\usepackage{algorithm}
\usepackage{algpseudocode}
\usepackage{booktabs}
\usepackage{enumitem}
\usepackage{xspace}
\usepackage[normalem]{ulem}
\usepackage{xcolor}
\usepackage{array}

\usepackage{graphicx}

\usepackage[title]{appendix}
\newenvironment{APPENDICES}{\begin{appendices}}{\end{appendices}}

\theoremstyle{plain}

\theoremstyle{definition}

\theoremstyle{remark}

\usepackage{natbib}
\bibpunct[, ]{(}{)}{,}{a}{}{,}

\usepackage{hyperref}
\hypersetup{
  colorlinks=true,
  linkcolor=darkgreen,
  filecolor=darkgreen,
  citecolor=darkgreen,     
  urlcolor=cyan,  
}

\usepackage{xr}

\allowdisplaybreaks

\definecolor{wisconsin-red}{rgb}{0.6,0,0}
\definecolor{darkgreen}{rgb}{0.2,0.6,0.2}
\definecolor{maroon}{rgb}{0.5, 0.0, 0.0}
\definecolor{violet}{rgb}{0.75, 0.0, 1.0}
\definecolor{lightgray}{gray}{0.9}
\definecolor{navyblue}{rgb}{0.0, 0.0, 0.5}
\definecolor{darkmidnightblue}{rgb}{0.0, 0.2, 0.4}
\definecolor{midnightblue}{rgb}{0.0,0.4,0.85}
\definecolor{Gray}{gray}{0.75}

\newcommand{\setItemSep}{\setlength\itemsep}

\newcommand{\revb}[1]{\textcolor{black}{#1}}

\newcolumntype{C}[1]{>{\centering\arraybackslash}p{#1}}
\newcolumntype{P}[1]{>{\raggedright\arraybackslash}p{#1}}
\newcolumntype{L}[1]{>{\raggedleft\arraybackslash}p{#1}}

\newcommand{\CombinedNeurADP}{\texttt{Combined-NeurADP}\xspace}
\newcommand{\VehicleNeurADP}{\texttt{Vehicle-NeurADP}\xspace}
\newcommand{\CentreNeurADP}{\texttt{Centre-NeurADP}\xspace}
\newcommand{\ILPThroughput}{\texttt{ILP-Throughput}\xspace}
\newcommand{\ILPVolume}{\texttt{ILP-Volume}\xspace}

\newcommand{\Time}{\mathcal{T}}
\newcommand{\VehiclesSet}{\mathcal{A}}
\newcommand{\CentresSet}{\mathcal{C}}
\newcommand{\SamplesSet}{\mathcal{B}}
\newcommand{\FeasibleSet}{\mathcal{F}}
\newcommand{\FeasibleDecisions}{\mathbf{X}}
\newcommand{\FeasibleDispatch}{\FeasibleSet^{\texttt{disp}}}
\newcommand{\FeasibleIdle}{\FeasibleSet^{\texttt{idle}}}
\newcommand{\FeasibleSkip}{\FeasibleSet^{\texttt{skip}}}
\newcommand{\AvailableVehicles}{\VehiclesSet^{\texttt{avail}}}

\newcommand{\CurrentTime}{t}
\newcommand{\StateNote}{S}
\newcommand{\EpochLength}{\delta}
\newcommand{\Deadline}{\tau}
\newcommand{\Capacity}{Q}
\newcommand{\NumAgents}{n}
\newcommand{\LocsToVisit}{L}
\newcommand{\DeadPeriod}{\Delta}
\newcommand{\VehicleState}{S^{\texttt{Veh}}}
\newcommand{\CentreState}{S^{\texttt{Ctr}}}
\newcommand{\PostDecisionState}{S^{\texttt{Post}}}
\newcommand{\Reward}{R}
\newcommand{\ImmediateReward}{r}
\newcommand{\DiscountRate}{\gamma}
\newcommand{\ExogenousInformation}{W}

\newcommand{\StepSize}{\alpha}
\newcommand{\Days}{N}
\newcommand{\SingleDay}{n}

\newcommand{\NullAction}{\varnothing}
\newcommand{\NoVehicle}{{\IndividualVehicle_{0}}}
\newcommand{\DecisionsVector}{\mathbf{x}}

\newcommand{\IndividualVehicle}{\texttt{$a$}}
\newcommand{\IndividualCentre}{\texttt{$c$}}
\newcommand{\IndividualSample}{\texttt{$b$}}

\newcommand{\VehicleReturn}{\text{return}}
\newcommand{\CentreVolume}{\text{volume}}
\newcommand{\CentreDeadline}{\text{deadline}}

\newcommand{\CentreSamples}{\text{samples}}

\newcommand{\SampleOrigin}{\text{origin}}
\newcommand{\SampleVolume}{\text{volume}}
\newcommand{\SampleDeadline}{\text{deadline}}

\newcommand{\LocationsSet}{\mathcal{L}}
\newcommand{\LocationIdx}{\ell}
\newcommand{\CentreSubset}{C'}
\newcommand{\PermSet}{\mathcal{P}}
\newcommand{\Route}{\sigma}
\newcommand{\OptimalRoute}{\sigma^*}
\newcommand{\UsedSamples}{\SamplesSet^{\texttt{used}}}
\newcommand{\CollectedVolume}{q}

\newcommand{\TravelTime}{\texttt{time}}
\newcommand{\RouteTime}{\texttt{RouteTime}}

\newcommand{\FirstTransition}{\texttt{statepost}}
\newcommand{\SecondTransition}{\texttt{statenext}}
\newcommand{\UsedSamplesAt}{\mathcal{U}}
\newcommand{\ExpiredSamplesAt}{\mathcal{E}}

\newcommand{\AuxFeatures}{\boldsymbol{\xi}}
\newcommand{\PostValueVeh}{\hat{V}^{\texttt{Veh}}}
\newcommand{\PostValueCtr}{\hat{V}^{\texttt{Ctr}}}
\newcommand{\TargetPostValueVeh}{\hat{V}^{\texttt{Veh-tgt}}}
\newcommand{\TargetPostValueCtr}{\hat{V}^{\texttt{Ctr-tgt}}}
\newcommand{\ActionScore}{\texttt{score}}

\title{\textbf{Dynamic Dispatching for Time-Sensitive Blood Sample Collection and Delivery}}

\author[1]{Arash Dehghan}
\author[1]{Aliaa Alnaggar}
\author[1]{Mucahit Cevik\thanks{Corresponding author: \href{mailto:mcevik@torontomu.ca}{mcevik@torontomu.ca}}}
\author[2]{Merve Bodur}

\affil[1]{Department of Mechanical, Industrial and Mechatronics Engineering, Toronto Metropolitan University, Toronto, Ontario, Canada}
\affil[2]{School of Mathematics and Maxwell Institute for Mathematical Sciences, University of Edinburgh, Edinburgh, UK} 

\date{}

\begin{document}\sloppy

\maketitle

\begin{abstract}
Hospitals and diagnostic laboratories rely on couriers to collect blood samples from geographically dispersed collection centres and deliver them for analysis before their short viability windows expire; late deliveries force costly re-collection and can delay diagnosis. We study the real-time dispatching of such a courier fleet, in which sample requests arrive stochastically at the centres throughout the day and a central dispatcher must repeatedly decide which vehicles to send, which centres each should visit, and whether to collect urgent samples immediately or consolidate them into later trips, subject to hard delivery deadlines and vehicle capacity limits.
Unlike static planning models, which fix routes before demand is known, and reactive heuristics, which respond only to the current backlog, our approach anticipates future arrivals when weighing immediate collection against consolidation. We formulate the problem as a Markov decision process and develop a neural approximate dynamic programming framework for centralized dispatch. The method introduces a dual value function decomposition that separately represents vehicle states and collection-centre states through neural networks trained on post-decision states. These learned estimates are integrated through a matching formulation that selects dispatch actions while balancing supply availability, demand urgency, and downstream opportunity cost. Computational experiments on a realistic Greater Toronto Area network compare the proposed policy with myopic baselines and ablation variants. Results show that the proposed dual value function policy raises the share of sample volume delivered on time by 1 to 9 percentage points over myopic baselines, with the largest gains under tight fleet, capacity, deadline, and routing constraints; these on-time gains, in turn, reduce reliance on costly external couriers.
\end{abstract}


\bigskip

\section{Introduction}

The global healthcare courier services market, valued at approximately USD 45.5 billion in 2025, is projected to reach around USD 63.7 billion by 2031, a compound annual growth rate of about 5.7\% that reflects the growing scale of medical logistics operations \citep{mordor_healthcare_courier_2025}. A central component of this infrastructure is the transportation of clinical specimens: on a daily basis, millions of biological samples are delivered from dispersed hospitals and clinics to centralized laboratories for testing \citep{wang2015routing}, an activity that continues to attract major investment. DHL Group, for instance, has committed 2 billion Euros through 2030 to expand its health logistics division~\citep{dhl_health_logistics_investment}.

The timely delivery of these specimens is critical: blood samples are highly perishable, with viability ranging from several hours to over 24 hours depending on the type of test \citep{anaya2021iterated}. However, diagnostic accuracy degrades well before this outer limit, as many hematological and biochemical analytes exhibit significant degradation within 12--24 hours depending on the test \citep{wu2017long}. In addition, delays in transporting samples for glucose testing can lead to falsely low readings due to ongoing cellular metabolism within the specimen \citep{mrazek2020errors}. To mitigate these effects, healthcare systems enforce much tighter operational deadlines on the maximum time between sample collection and laboratory arrival---commonly 90 minutes for standard blood work \citep{zufferey2016dynamic} and no more than two hours for broader clinical analyses \citep{grasas2014improvement}.

Within this context, failures in the pre-analytical phase of laboratory testing are well documented. Pre-analytical errors, encompassing specimen collection, handling, transportation, and storage, account for approximately 60--70\% of all problems in laboratory diagnostics \citep{lippi2011preanalytical}. While these errors arise from multiple sources, transportation delays represent a critical and operationally controllable component, directly influenced by courier routing and dispatching decisions. In a separate analysis of reported laboratory incidents, 60\% were associated with a potential diagnostic error, most commonly involving delays in the diagnostic process (50.5\%), along with misdiagnoses (32.1\%) and missed diagnoses (17.3\%) \citep{van2023nature}. Moreover, these failures can extend beyond diagnostic inefficiencies to adverse patient outcomes. In an analysis of 684 laboratory-related incidents across multiple healthcare facilities, 8\% resulted in additional treatment requirements or temporary or permanent harm to the patient \citep{tran2020pre}.

The financial burden is also substantial. Industry analyses estimate that in a typical three- to four-hospital system, specimen mishandling and re-collection costs can exceed 1 million USD annually, with per-incident costs ranging from \$350 to \$5,000 depending on the specimen type \citep{MedSpeed_2025}. Indeed, lost or delayed delivery of specimens has been identified as the most common problem jeopardizing patient safety in clinical laboratory operations \citep{astion2003classifying}. Transporting clinical specimens under these constraints represents a significant logistical burden to healthcare systems, particularly for cases requiring prompt courier services \citep{wang2015routing}.

Given these clinical and economic stakes, the operations research literature has devoted considerable attention to routing vehicles that collect biomedical samples and deliver them to a central laboratory before their viability expires. This body of work is predominantly \textit{static}: requests, time windows, and volumes are assumed known in advance, so a complete collection plan is optimized offline before operations begin \citep{grasas2014improvement, toschi2018fix, anaya2021iterated, benini2022mathematical, mazzanti2024iterative}.Only \mbox{\citet{wang2015routing}} and \mbox{\citet{zufferey2016dynamic}} incorporate elements of uncertainty. The former models uncertainty through a two-stage stochastic program in which all demand scenarios are resolved prior to vehicle departure, while the latter employs a greedy insertion heuristic that reactively accommodates dynamic requests as they arise during the day.
To our knowledge, no prior work addresses the fully dynamic, stochastic variant of the problem using non-myopic solution approaches, despite the substantial gains such methods have demonstrated in other transportation domains~\citep{shah2020neural, dehghan2025dynamic, dehghan2026neural}.

In this paper, we study real-time dispatching for blood sample collections under stochastic demand. Sample requests arrive dynamically at collection centres throughout the day, and a central laboratory must make sequential dispatch decisions under uncertainty, determining which vehicles to deploy, which centres to visit, and whether to consolidate pickups or dispatch immediately for urgent samples. This setting is closely related to the multi-vehicle routing problem studied by \citet{zufferey2016dynamic}, where a cost-effective primary fleet performs most collections while unmet demand is handled through costlier alternatives. However, existing approaches rely on myopic, reactive heuristics that do not anticipate future demand or account for the downstream impact of current decisions. To address this limitation, we formulate the dispatching problem as a Markov decision process (MDP) and develop a Neural Approximate Dynamic Programming (NeurADP) framework for real-time vehicle dispatch. Our model captures key operational features of specimen collection systems, including stochastic arrivals, hard perishability deadlines, multi-stop routing with capacity constraints, and centralized real-time decision-making. We introduce a dual value function decomposition that maintains separate neural networks for vehicle and collection-centre states, enabling independent modeling of supply and demand dynamics. These components are integrated through a matching formulation. We evaluate our approach on a realistic network based on the Greater Toronto Area and compare it against multiple benchmark policies, including myopic baselines and ablation variants.

The contributions of our work are summarized as follows:
\begin{itemize}[leftmargin=1.2em]
\item We formulate a novel MDP model for the blood sample collection and delivery problem that accounts for uncertainty in sample arrivals within a centralized collection system with perishability deadlines and multi-stop vehicle routing. To our knowledge, this is the first study to address the fully dynamic variant of this problem with a non-myopic, learning-based solution approach.

\item We propose a dual value function architecture for NeurADP that decomposes system-level value into vehicle-specific and centre-specific components, each approximated by a separate neural network trained on entity-level post-decision states. This decomposition extends single-value-function approaches by allowing the policy to learn supply-side and demand-side dynamics separately and combine them through a matching formulation.

\item We conduct an extensive computational study showing that the proposed approach improves volume served by 1\% to 9\% relative to myopic baselines across a range of operational settings, with the largest gains under tight constraints. Beyond the headline service rate, we trace the advantage to how the policy balances route consolidation against timely collection---sustaining frequent, efficient dispatches while keeping spatial coverage across centres balanced---and note that the resulting on-time gains reduce reliance on costly external couriers. An ablation study confirms that both value function components contribute significantly and that their combination outperforms either component alone in all but the simplest routing setting. 

\item We provide managerial insights for specimen collection operations. Specifically, we analyze the impact of fleet size, vehicle capacity, sample perishability deadlines, and route complexity, and show that the proposed policy performs well by adaptively balancing route consolidation with timely collection.
\end{itemize}

The rest of the paper is structured as follows. Section~\ref{literaturereview} reviews the most relevant prior work and situates our study within the existing literature. Section~\ref{problemdescription} describes the blood sample delivery problem. Section~\ref{modelformulation} introduces the formal MDP formulation. Section~\ref{solutionmethodology} details the proposed NeurADP solution methodology. 
Section~\ref{experimentalsetup} describes the experimental design, neural network architectures, and benchmark policies. Section~\ref{results} reports and analyzes the outcomes of our computational study, and Section~\ref{conclusion} summarizes the main contributions and outlines directions for future work.

\section{Literature Review} \label{literaturereview}

We review three streams of work that frame our contribution: blood and biomedical sample transportation, which establishes the application domain; Approximate Dynamic Programming (ADP) for dynamic dispatching and delivery, which provides the methodological foundation; and value function decomposition for multi-entity systems, which motivates the dual vehicle-and-centre architecture that constitutes our main methodological contribution.

\subsection{Blood and Biomedical Sample Transportation}
Table~\ref{table:references} shows the most relevant studies along key problem and methodological dimensions. Among the columns, ``Dispatch Optimized'' indicates whether the method explicitly optimizes vehicle-to-location matching at each epoch, and ``Non-Myopic'' indicates policy anticipating future arrivals rather than optimizing only the immediate reward; the remaining columns are self-explanatory.

\setlength{\tabcolsep}{2.5pt}
\renewcommand{\arraystretch}{0.905}
\begin{table}[!ht]
\centering
\caption{Comparison of blood and biomedical sample transportation literature.}\label{table:references}
\resizebox{0.89\textwidth}{!}{
\begin{tabular}{P{0.32\textwidth}C{0.11\textwidth}C{0.11\textwidth}C{0.11\textwidth}C{0.11\textwidth}C{0.13\textwidth}C{0.11\textwidth}C{0.11\textwidth}}
\toprule
\textbf{Study} & \textbf{Stochastic Arrivals} & \textbf{Online Dispatch} & \textbf{Multi-Trip} & \textbf{Deadlines} & \textbf{Dispatch Optimized} & \textbf{Non-Myopic} & \textbf{ADP / VFA} \\
\midrule
\citet{grasas2014improvement}        &            &            &            & \checkmark &            &            &            \\
\citet{wang2015routing}              & \checkmark &            & \checkmark & \checkmark &            &            &            \\
\citet{zufferey2016dynamic}          & \checkmark & \checkmark & \checkmark & \checkmark &            &            &            \\
\citet{toschi2018fix}     &            &            & \checkmark & \checkmark &            &            &            \\
\citet{anaya2021iterated}      &            &            & \checkmark & \checkmark &            &            &            \\
\citet{benini2022mathematical}       &            &            & \checkmark & \checkmark &            &            &            \\
\citet{mazzanti2024iterative}        &            &            & \checkmark & \checkmark &            &            &            \\
\citet{ocampo2025iterative}   &            &            & \checkmark & \checkmark &            &            &            \\
\midrule
\textbf{Our Work}                    & \checkmark & \checkmark & \checkmark & \checkmark & \checkmark & \checkmark & \checkmark \\
\bottomrule
\end{tabular}
}
\end{table}

Early research on blood and biomedical sample transportation has focused predominantly on static, deterministic routing, in which all requests, time windows, and volumes are assumed to be known in advance and routes are constructed and fixed before operations. \citet{grasas2014improvement} model blood sample collection for laboratories in Catalonia, Spain, as a capacitated, time-constrained open vehicle routing problem solved with a biased random-key genetic algorithm. A closely related stream developed with the Quebec Ministry of Health progressively enriches this line of research. \citet{toschi2018fix} propose a fix-and-optimize variable neighbourhood search that treats the number of visits per centre as a decision; \citet{anaya2021iterated} add interdependent pickups and flexible center opening hours; \citet{mazzanti2024iterative} introduce a rolling-horizon matheuristic for large-scale realistic instances; and \citet{ocampo2025iterative} recast the problem as service network design over a time-expanded network, yielding the first exact approach at scale. In a related setting, \citet{benini2022mathematical} incorporate sample stabilization at spoke centres, inter-vehicle transfers, and multi-depot routing. Despite their sophistication, all assume that sample generation times are known in advance and therefore cannot accommodate dynamically arriving requests.

A smaller body of work incorporates demand uncertainty but does not account for real-time anticipation. \citet{wang2015routing} address healthcare courier delivery with stochastic urgent demand through a two-stage stochastic program that hedges master routes against demand scenarios; however, all uncertainty is resolved before courier departure and the recourse stage relies on greedy insertion, leaving dispatching effectively offline. \citet{zufferey2016dynamic} more closely approximates real-time operation through the Medical Vehicle Routing Problem, motivated by a Geneva hospital laboratory and incorporating a heterogeneous car-and-scooter fleet, time-dependent travel times, laboratory-based deadlines, and dynamically arriving requests; their greedy-insertion and local-search heuristic diverts vehicles as requests appear but remains purely reactive, acting on the current state without anticipating future arrivals. At a strategic level, \citet{mousavi2021designing} study multi-objective blood supply chain network design under uncertainty using metaheuristics and chance constraints, which complements rather than directly addresses operational dispatching. In short, we are not aware of prior work that tackles the fully dynamic, stochastic dispatching problem with a non-myopic policy. This gap motivates ADP, which has delivered strong anticipatory, value-based policies in related transportation domains and which we review next.

\subsection{Approximate Dynamic Programming for Dynamic Dispatching and Delivery}

Early ADP approaches to dynamic dispatching typically represent the post-decision value function using either lookup tables or state-aggregation schemes, thereby combining similar states to limit the size of the approximation. Such approaches have been applied to long-haul trucking \citep{simao2010approximate}, ambulance dispatching and relocation \citep{schmid2012solving}, electric-vehicle fleet management \citep{al2020approximate}, last-mile ride-sharing \citep{agussurja2019state}, and bike-sharing inventory routing \citep{brinkmann2019dynamic}. Subsequent work extends ADP to richer decisions and larger systems. In dynamic dial-a-ride, \citet{heitmann2023combining} combine a value-function approximation for service-offering decisions with scenario-based routing, while \citet{heitmann2024accelerating} improve convergence by progressively moving from lower- to higher-dimensional state representations. Most relevant to our setting, \citet{van2019delivery} study same-day delivery dispatching for urban consolidation centres and embed a linear parametric approximation of the value function, constructed from manually specified basis functions, within an integer program. Their consolidation-versus-urgency trade-off is directly analogous to the one considered here. These methods nevertheless depend on manually designed state aggregations or feature representations, which require substantial domain expertise and may become difficult to scale as the number of entities and state attributes increases.

This limitation motivated neural VFA, which learns state representations directly from data. \citet{shah2020neural} introduced Neural ADP (NeurADP) for ride-pooling problem, approximating post-decision state values with neural networks and making dispatch decisions through an integer program. The framework has since been extended to joint matching and pricing in crowd-shipping \citep{dehghan2026joint} and to flexible pickup and drop-off locations \citep{jiang2026optimizing}, refined with value decomposition strategies that better capture inter-vehicle dependencies \citep{bose2021conditional, hao2022hierarchical}, which we revisit in Section~\ref{sec:vfdecomp}. 
A common structural limitation nonetheless persists: these approaches decompose the value function along a single entity type, typically vehicles, so demand-side dynamics can be represented only indirectly---as input features to the vehicle value function rather than through value terms of their own. In blood sample collection, system performance hinges on the interaction between two fundamentally different entities, vehicles and collection centres, with distinct states and temporal dynamics, a structure that a single-entity decomposition cannot capture adequately.

\subsection{Value Function Decomposition for Multi-Entity Systems} \label{sec:vfdecomp}

Value function decomposition has been pursued along several dimensions. For
multiattribute resource management, \citet{george2008value} combine value
estimates obtained from progressively coarser representations of a resource's
attribute vector, balancing aggregation bias against sampling error. In
electric-vehicle fleet management, \citet{al2020approximate} similarly estimate
vehicle-state values using a hierarchy of coarsened attribute representations,
where a vehicle is characterized by its availability time, location, and battery
level. Both approaches reduce dimensionality within the state representation of
a single resource class, rather than assigning separate value components to
different types of system entities. By contrast, NeurADP approach proposed by \citet{shah2020neural} decomposes the joint post-decision value into a sum of
vehicle-specific value functions, each evaluated using the focal vehicle's
post-decision state and the other vehicles' states as context. This per-vehicle decomposition has been refined in ride-pooling: \citet{bose2021conditional} correct for inter-vehicle competition using the conditional action probabilities of neighbouring agents, and \citet{hao2022hierarchical} introduce a hierarchical mixing architecture that combines per-vehicle values within and across spatial clusters; alternatively, \citet{yu2019integrated} split ride-pooling into separate matching and routing subproblems. All of these refinements operate along a single dimension, the vehicle fleet, because their demand-side entities (passenger requests) are transient, leaving no persistent demand state to track across epochs. In blood sample collection, by contrast, collection centres are persistent: samples accumulate with rising urgency, creating demand-side dynamics that cannot be captured from the vehicle side alone.

Among the prior studies, the meal-delivery problem of \citet{neria2025restaurant} is most relevant to our multi-entity setting as it involves two distinct resource types---cooks who prepare orders and vehicles that deliver them. Their approach, however, does not decompose the value function by resource type. Instead, system-wide summary statistics describing both cooks and vehicles are combined into a single fixed-dimensional feature vector, which is then passed to one neural network. This representation is designed to facilitate transfer across instances of different sizes. Because the information is aggregated at the system level, however, it does not isolate the marginal value of an individual resource and therefore does not provide the entity-level value estimates required by a structured matching formulation. In contrast, our NeurADP framework for blood sample collection maintains separate value functions for vehicles and collection centres, each approximated by a dedicated neural network. The vehicle component captures supply-side dynamics, whereas the centre component captures demand-side dynamics. The entity-level value estimates are then integrated within the matching formulation to evaluate joint dispatch and collection decisions. 

\section{Problem Description} \label{problemdescription}

We consider the dynamic collection of perishable blood samples from geographically distributed collection centres and their delivery to a central hospital laboratory for analysis. Operations take place over a single business day, partitioned into discrete decision epochs of length~$\EpochLength$ minutes. At each epoch, a central dispatcher determines which available vehicles to dispatch, the collection centres assigned to each dispatched vehicle, and the sequence in which those centres are visited. These decisions must balance the urgency of samples approaching their delivery deadlines against the operational efficiencies obtained by consolidating pickups into multi-centre routes.

More specifically, the collection network consists of a central hospital and a set $\CentresSet = \{1,2,\ldots,|\CentresSet|\}$ of collection centres. The hospital serves as both the depot from which vehicles depart and the destination to which all collected samples are delivered. A fleet of $\NumAgents$ identical vehicles is stationed at the hospital, each with a maximum carrying capacity of $\Capacity$ sample-volume units. The homogeneous-fleet assumption reflects settings in which a standardized urban courier mode, such as scooters, mopeds, or bicycles, is used for time-sensitive medical sample transportation. 
Vehicles operate in a multi-trip fashion: after completing a collection route and returning to the hospital, a vehicle becomes available for reassignment in a subsequent decision epoch.

Blood samples arrive stochastically at collection centres throughout the day. The operating day consists of an active arrival period followed by a terminal clearance period of $\DeadPeriod$ minutes during which no new samples arrive but vehicles continue to operate and collect outstanding samples, a structure commonly used in dynamic specimen collection settings \citep{zufferey2016dynamic}. During the active period, the number of samples arriving at each epoch follows a non-homogeneous Poisson process, consistent with prior work modeling time-varying demand in healthcare logistics \citep{wang2015routing}. The spatial distribution of arrivals across centres varies stochastically over time to reflect heterogeneous and shifting patient demand. Each sample~$\IndividualSample \in \SamplesSet$ is characterized by its origin centre, its arrival time, and its volume. Upon arrival, each sample is assigned a hard delivery deadline of $\Deadline$ minutes from its arrival time, reflecting the time-sensitive nature of clinical specimen viability and commonly enforced operational constraints in laboratory logistics \citep{zufferey2016dynamic, benini2022mathematical, ocampo2025iterative}. If a sample is not collected and delivered to the hospital before its deadline expires, it exits the primary-fleet system and is counted as lost for the purposes of evaluating fleet performance; operationally, such samples are assumed to be served by an external courier at higher cost, consistent with the two-tier service models studied in prior work \citep{zufferey2016dynamic}. Any samples remaining at collection centres at the end of the operating day are also counted as lost.

At each decision epoch~$\CurrentTime$, the dispatcher observes the current system state, including the status of each vehicle and the outstanding sample backlog at each collection centre, with all arrivals realized by the start of the epoch. Based on this information, the dispatcher selects a set of vehicle dispatch decisions. Each dispatch corresponds to a collection route: a vehicle departs the hospital, visits a subset of collection centres, collects available samples subject to its remaining capacity, and returns to the hospital. The visiting order of centres is chosen to minimize total travel time over the selected centres. A fixed service time is incurred at each visited centre to account for sample handoff and loading. Upon returning to the hospital, the vehicle deposits all collected samples and becomes available for dispatch in a future epoch. 
Each vehicle may visit at most $\LocsToVisit$ collection centres in a single trip, reflecting route-duration and operational-complexity limits commonly imposed in healthcare logistics \citep{wang2015routing, zufferey2016dynamic}. By preventing excessively long multi-stop routes, this restriction helps preserve vehicle availability for newly arriving or increasingly urgent samples in subsequent epochs and reduces the risk of missed deadlines and reliance on external courier services.

When a vehicle visits a centre, it prioritizes the most urgent samples, collecting them in order of earliest deadline, which reflects standard handling practices for perishable medical specimens. Any samples that cannot be accommodated due to capacity limitations remain at the centre for future collection provided their deadlines have not expired. 
At each decision epoch, centres that are not visited retain their samples, which continue to age toward their deadlines, and vehicles that are not dispatched remain idle at the hospital. This allows the dispatcher to strategically delay dispatch decisions when future consolidation opportunities are expected to improve efficiency, while still managing the risk of sample expiration.

\section{Model Formulation} \label{modelformulation}

This section formalizes the blood sample collection and delivery problem described in Section~\ref{problemdescription} as an MDP. We first define the state and decision variables, then specify the reward structure, exogenous information, transition function, and finite-horizon objective.

\subsection{State Variables}

At each decision epoch $\CurrentTime \in \Time$, the system state is given by the tuple $\StateNote_\CurrentTime = (\VehicleState_\CurrentTime, \CentreState_\CurrentTime)$, where $\VehicleState_\CurrentTime$ captures the state of the fleet and $\CentreState_\CurrentTime$ captures the state of the collection centres, including the blood samples currently buffered at each. We identify each collection centre with its index in $\{1, \ldots, |\CentresSet|\}$, so that $\CentresSet$ is a set of location indices, and let $\LocationsSet = \{0\} \cup \CentresSet$ denote the set of physical locations, where $\LocationIdx = 0$ is the central hospital. The travel time from location $\LocationIdx$ to location $\LocationIdx'$ is denoted $\TravelTime(\LocationIdx, \LocationIdx')$ and includes a fixed service time incurred upon arrival at any collection centre but not at the hospital. We represent each vehicle $\IndividualVehicle \in \VehiclesSet$ by a one-dimensional attribute $\IndividualVehicle = (\IndividualVehicle_\VehicleReturn)$, where $\IndividualVehicle_\VehicleReturn \in \mathbb{R}_{\geq 0}$ denotes the duration (in minutes) before the vehicle completes its current route and returns to the hospital. A vehicle with $\IndividualVehicle_\VehicleReturn = 0$ is idle and eligible for dispatch at the current epoch. The fleet state $\VehicleState_\CurrentTime$ collects the attribute $\IndividualVehicle_\VehicleReturn$ of every vehicle $\IndividualVehicle \in \VehiclesSet$.

We represent each blood sample $\IndividualSample \in \SamplesSet_\CurrentTime$ by a three-dimensional attribute vector $\IndividualSample = (\IndividualSample_\SampleOrigin, \IndividualSample_\SampleVolume, $ $\IndividualSample_\SampleDeadline)$, where $\IndividualSample_\SampleOrigin \in \CentresSet$ is the origin centre, $\IndividualSample_\SampleVolume \in \mathbb{R}_{>0}$ is the sample volume, and $\IndividualSample_\SampleDeadline \in \mathbb{R}_{>0}$ is the hard delivery deadline. Here, $\SamplesSet_\CurrentTime$ denotes the set of all blood samples currently in the system at epoch $\CurrentTime$. When a sample arrives at epoch $\CurrentTime$, its deadline is initialized as
\begin{equation}
\IndividualSample_\SampleDeadline = \CurrentTime + \Deadline,
\label{eq:deadline-assignment}
\end{equation}
and the sample is placed at its origin centre. A sample is removed from the system and counted as lost whenever it can no longer feasibly be collected and delivered to the hospital before its deadline, or if it remains at a collection centre at the end of the operating day.

We represent the state of each collection centre $\IndividualCentre \in \CentresSet$ by the three-dimensional attribute vector $(\IndividualCentre_\CentreSamples, \IndividualCentre_\CentreVolume, \IndividualCentre_\CentreDeadline)$, where $\IndividualCentre_\CentreSamples \subseteq \SamplesSet_\CurrentTime$ is the set of samples currently buffered at the centre and the remaining two attributes are summary features derived from it: $\IndividualCentre_\CentreVolume = \sum_{\IndividualSample \in \IndividualCentre_\CentreSamples} \IndividualSample_\SampleVolume$ is the total stored volume, and $\IndividualCentre_\CentreDeadline = \min_{\IndividualSample \in \IndividualCentre_\CentreSamples} \IndividualSample_\SampleDeadline$ is the earliest deadline among buffered samples. For a centre without available samples, we set $\IndividualCentre_\CentreVolume = 0$ and $\IndividualCentre_\CentreDeadline = \infty$. The centre state $\CentreState_\CurrentTime$ collects the attribute tuple of every centre $\IndividualCentre \in \CentresSet$. By construction, the buffered sample sets partition the outstanding samples, so $\SamplesSet_\CurrentTime = \bigcup_{\IndividualCentre \in \CentresSet} \IndividualCentre_\CentreSamples$.

\subsection{Decision Variables}
At each decision epoch $\CurrentTime \in \Time$, the dispatcher selects a joint assignment that simultaneously determines (i) which available vehicles to dispatch, (ii) which collection centres each dispatched vehicle visits, and (iii) how collection along the resulting route is executed given the vehicle's capacity and the buffered samples' deadlines. We formalize this joint decision as the selection of exactly one elementary action per vehicle and exactly one coverage action per centre from a pool of feasible alternatives, subject to mutual exclusivity constraints that avoid double-counting. An elementary vehicle action is a pair $(\IndividualVehicle, \CentreSubset)$, where \revb{$\IndividualVehicle \in \AvailableVehicles_\CurrentTime := \{\IndividualVehicle \in \VehiclesSet : \IndividualVehicle_\VehicleReturn = 0\}$ is a vehicle available for dispatch at epoch $\CurrentTime$} and $\CentreSubset \subseteq \CentresSet$ is a non-empty subset of at most $\LocsToVisit$ collection centres. For each such pair, we determine the order in which the centres are visited by solving a shortest-path problem over all permutations of $\CentreSubset$:
\begin{equation}
\OptimalRoute_{\IndividualVehicle, \CentreSubset} \in \arg\min_{\Route \in \PermSet(\CentreSubset)} \left\{ \TravelTime(0, \Route_1) + \sum_{j=1}^{|\CentreSubset|-1} \TravelTime(\Route_j, \Route_{j+1}) + \TravelTime(\Route_{|\CentreSubset|}, 0) \right\},
\label{eq:tsp-route}
\end{equation}
where $\PermSet(\CentreSubset)$ denotes the set of all permutations of the centres in $\CentreSubset$, the route begins and ends at the hospital ($\LocationIdx = 0$), and $\Route_j$ is the $j^{\text{th}}$ centre visited in the ordering. Let $\RouteTime(\IndividualVehicle, \CentreSubset)$ denote the total travel time along this optimal route.

Given a vehicle action $(\IndividualVehicle, \CentreSubset)$, samples are collected greedily in earliest-deadline-first order among the buffered samples available at the visited centres, and we write $\UsedSamples_{\IndividualVehicle, \CentreSubset} \subseteq \bigcup_{\IndividualCentre \in \CentreSubset} \IndividualCentre_\CentreSamples$ for the resulting set of samples collected along the route. This selection is subject to two feasibility conditions. First, the total volume of samples collected must not exceed the vehicle's capacity:
\begin{equation}
\CollectedVolume_{\IndividualVehicle, \CentreSubset} := \sum_{\IndividualSample \in \UsedSamples_{\IndividualVehicle, \CentreSubset}} \IndividualSample_\SampleVolume \leq \Capacity.
\label{eq:capacity-constraint}
\end{equation}
Second, every collected sample must be deliverable to the hospital before its deadline expires:
\begin{equation}
\CurrentTime + \RouteTime(\IndividualVehicle, \CentreSubset) \leq \IndividualSample_\SampleDeadline, \qquad \forall \IndividualSample \in \UsedSamples_{\IndividualVehicle, \CentreSubset}.
\label{eq:deadline-constraint}
\end{equation}

In addition to vehicle actions, the decision includes two classes of \revb{null} actions, which ensure that every vehicle and every centre is accounted for in the assignment. For each vehicle $\IndividualVehicle \in \VehiclesSet$, the vehicle-idle action $(\IndividualVehicle, \NullAction)$ represents not dispatching the vehicle at the current epoch: an available vehicle remains idle at the hospital, while a vehicle still en route continues its current trip with its time until return decrementing by one epoch length. For each centre \revb{$\IndividualCentre \in \CentresSet$}, the centre-skip action \revb{$(\NoVehicle, \{\IndividualCentre\})$, where $\NoVehicle \notin \VehiclesSet$ denotes a dummy vehicle,} represents leaving the centre unattended at the current epoch, with its stored samples ageing toward their deadlines.

The elementary actions available at epoch $\CurrentTime$ therefore fall into three families,
\begin{align*}
\FeasibleDispatch_\CurrentTime(\StateNote_\CurrentTime) &= \big\{ (\IndividualVehicle, \CentreSubset) \;:\; \IndividualVehicle \in \AvailableVehicles_\CurrentTime, \; \NullAction \neq \CentreSubset \subseteq \CentresSet, \; |\CentreSubset| \leq \LocsToVisit, \; \eqref{eq:capacity-constraint}\text{--}\eqref{eq:deadline-constraint} \text{ hold} \big\}, \\
\FeasibleIdle &= \big\{ (\IndividualVehicle, \NullAction) \;:\; \IndividualVehicle \in \VehiclesSet \big\}, \\
\FeasibleSkip &= \big\{ (\NoVehicle, \{\IndividualCentre\}) \;:\; \IndividualCentre \in \CentresSet \big\},
\end{align*}
and we write $\FeasibleSet_\CurrentTime(\StateNote_\CurrentTime) = \FeasibleDispatch_\CurrentTime(\StateNote_\CurrentTime) \cup \FeasibleIdle \cup \FeasibleSkip$ for their union, which is disjoint because the three families are distinguished by whether the centre subset is empty and whether a vehicle is present. Dispatch actions are generated only for vehicles that are available at the current epoch, whereas idle actions are defined for every vehicle, including those still en route, and skip actions for every centre; consequently every vehicle and every centre can always be covered.

Let $x_{\CurrentTime, (\IndividualVehicle, \CentreSubset)} \in \{0, 1\}$ indicate whether the elementary action $(\IndividualVehicle, \CentreSubset) \in \FeasibleSet_\CurrentTime(\StateNote_\CurrentTime)$ is selected at epoch $\CurrentTime$, and let $\DecisionsVector_\CurrentTime = \big(x_{\CurrentTime, (\IndividualVehicle, \CentreSubset)}\big)_{(\IndividualVehicle, \CentreSubset) \in \FeasibleSet_\CurrentTime(\StateNote_\CurrentTime)}$ denote the full decision vector. The decision must satisfy
\begin{subequations}
\label{eq:assignment-constraints}
\begin{align}
& \sum_{\CentreSubset \,:\, (\IndividualVehicle, \CentreSubset) \in \FeasibleSet_\CurrentTime(\StateNote_\CurrentTime)} x_{\CurrentTime, (\IndividualVehicle, \CentreSubset)} = 1, & \quad \forall \IndividualVehicle \in \VehiclesSet, \label{eq:vehicle-constraint} \\
& \sum_{(\IndividualVehicle, \CentreSubset) \in \FeasibleSet_\CurrentTime(\StateNote_\CurrentTime) \,:\, \IndividualCentre \in \CentreSubset} x_{\CurrentTime, (\IndividualVehicle, \CentreSubset)} = 1, & \quad \forall \IndividualCentre \in \CentresSet. \label{eq:centre-constraint}
\end{align}
\end{subequations}
Constraints~\eqref{eq:vehicle-constraint} requires that each vehicle be covered by exactly one elementary action, which is either a dispatch to a centre subset or an idle action. Constraints~\eqref{eq:centre-constraint} requires that each collection centre be covered by exactly one elementary action, which is either through a dispatched vehicle whose route includes it or through a centre-skip action. We denote by $\FeasibleDecisions_\CurrentTime(\StateNote_\CurrentTime)$ the set of binary decision vectors $\DecisionsVector_\CurrentTime$ satisfying~\eqref{eq:assignment-constraints}, which constitutes the feasible joint assignments at epoch $\CurrentTime$.

\revb{Because every vehicle and every centre is covered exactly once, any objective that is additively separable across elementary actions receives exactly one term per vehicle and one per centre; Section~\ref{sec:dualvfa} exploits this to evaluate a joint assignment through per-entity value estimates.}

\subsection{Rewards}

The immediate reward at each decision epoch is the total volume of blood samples selected for collection in the current epoch, reflecting the primary operational objective of maximizing the total sample volume delivered to the hospital within their deadlines. \revb{For a dispatch action $(\IndividualVehicle, \CentreSubset) \in \FeasibleDispatch_\CurrentTime(\StateNote_\CurrentTime)$,} the action reward is
\begin{equation}
\ImmediateReward(\IndividualVehicle, \CentreSubset) = \CollectedVolume_{\IndividualVehicle, \CentreSubset},
\label{eq:action-reward}
\end{equation}
namely, the total volume of samples collected along the optimal route $\OptimalRoute_{\IndividualVehicle, \CentreSubset}$. Idle and skip actions contribute no reward, i.e., \revb{$\ImmediateReward(\IndividualVehicle, \NullAction) = \ImmediateReward(\NoVehicle, \{\IndividualCentre\}) = 0$}. The total immediate reward obtained at epoch $\CurrentTime$ under decision vector $\DecisionsVector_\CurrentTime$ is
\begin{equation}
\Reward_\CurrentTime(\DecisionsVector_\CurrentTime) = \sum_{(\IndividualVehicle, \CentreSubset) \in \FeasibleSet_\CurrentTime(\StateNote_\CurrentTime)} x_{\CurrentTime, (\IndividualVehicle, \CentreSubset)} \cdot \ImmediateReward(\IndividualVehicle, \CentreSubset).
\label{eq:epoch-reward}
\end{equation}

\subsection{Random Information}
\label{sec:random_info}
At the end of each epoch $\CurrentTime$, the system observes the collection of new blood samples that have arrived at collection centres during the interval $(\CurrentTime, \CurrentTime + \EpochLength]$. This collection constitutes the exogenous information driving the system's evolution, which we denote by
\begin{equation}
\ExogenousInformation_{\CurrentTime + 1} = \{\IndividualSample_1, \IndividualSample_2, \ldots, \IndividualSample_{N_{\CurrentTime + 1}}\},
\label{eq:exogenous}
\end{equation}
where $N_{\CurrentTime + 1} \geq 0$ is the random number of samples arriving during the interval. Each arrival $\IndividualSample \in \ExogenousInformation_{\CurrentTime + 1}$ carries a random origin centre $\IndividualSample_\SampleOrigin \in \CentresSet$ and a random volume $\IndividualSample_\SampleVolume \in \mathbb{R}_{> 0}$, while its delivery deadline $\IndividualSample_\SampleDeadline$ is determined deterministically from the arrival epoch via~\eqref{eq:deadline-assignment}. Once observed, each arrival is placed at its origin centre and becomes part of that centre's buffered sample set $\IndividualCentre_\CentreSamples$, so that at epoch $\CurrentTime + 1$ the newly arrived samples are available for inclusion in any future dispatch decision. At the start of the operating horizon we assume $\ExogenousInformation_0 = \varnothing$, i.e., no samples have accumulated before decision-making begins. The specific arrival process---including the intensity of arrivals, the spatial distribution of origins, and the distribution of sample volumes---is specified in Section~\ref{sec:dataset}, where it is calibrated to empirical outpatient clinic data.

Because sample arrivals reflect patient flow to collection centres rather than any feedback from the dispatcher's actions, the distribution of $\ExogenousInformation_{\CurrentTime + 1}$ is exogenous conditional on the epoch: it may vary over time, but it does not depend on the current decision vector $\DecisionsVector_\CurrentTime$ or on past dispatch decisions. The mechanism by which $\ExogenousInformation_{\CurrentTime + 1}$ updates the system state to produce the next pre-decision state $\StateNote_{\CurrentTime + 1}$ is made explicit in the transition function presented next.

\subsection{Transition Function}

Following \citet{powell2007approximate}, we introduce a post-decision state and decompose the transition from $\StateNote_\CurrentTime$ to $\StateNote_{\CurrentTime + 1}$ into a deterministic transition from the pre-decision state to the post-decision state, induced by $\DecisionsVector_\CurrentTime$, and a stochastic transition from the post-decision state to the next pre-decision state, driven by the realization of exogenous information $\ExogenousInformation_{\CurrentTime + 1}$:
\begin{equation}
\PostDecisionState_\CurrentTime = \FirstTransition(\StateNote_\CurrentTime, \DecisionsVector_\CurrentTime), \qquad
\StateNote_{\CurrentTime + 1} = \SecondTransition(\PostDecisionState_\CurrentTime, \ExogenousInformation_{\CurrentTime + 1}).
\label{eq:transition-decomp}
\end{equation}
The post-decision state, $\PostDecisionState_\CurrentTime$, captures the system configuration immediately after the action is selected but before new samples arrive. Both state components evolve during this step: the return time of each vehicle advances to reflect its newly selected route or its continuing trip, and the sample set at each centre shrinks as samples are collected or expire.

For each vehicle $\IndividualVehicle \in \VehiclesSet$, time left until return post-decision is updated according to
\begin{equation}
\IndividualVehicle_\VehicleReturn^{\texttt{Post}} =
\begin{cases}
\max\!\left(\RouteTime(\IndividualVehicle, \CentreSubset) - \EpochLength,\, 0\right), & \text{if } \IndividualVehicle \text{ is dispatched along route } \OptimalRoute_{\IndividualVehicle, \CentreSubset}, \\
\max\!\left(\IndividualVehicle_\VehicleReturn - \EpochLength,\, 0\right), & \text{if } \IndividualVehicle \text{ is left idle or remains en route}.
\end{cases}
\label{eq:vehicle-transition}
\end{equation}
In both branches, one epoch length is deducted because the post-decision clock reflects the vehicle's state at the start of epoch $\CurrentTime + 1$, by which time one epoch has already elapsed: a dispatched vehicle has begun traversing its optimal route, and an idle or en-route vehicle has advanced one epoch closer to returning to the hospital.

For each collection centre $\IndividualCentre \in \CentresSet$, the post-decision sample set is obtained by removing the samples assigned to a selected route at the current epoch and any remaining samples whose deadlines can no longer be met by a future dispatch:
\begin{equation}
\IndividualCentre_\CentreSamples^{\texttt{Post}} = \IndividualCentre_\CentreSamples \setminus \left( \UsedSamplesAt_\CurrentTime(\IndividualCentre) \cup \ExpiredSamplesAt_\CurrentTime(\IndividualCentre) \right),
\label{eq:centre-post-transition}
\end{equation}
where $\UsedSamplesAt_\CurrentTime(\IndividualCentre) \subseteq \IndividualCentre_\CentreSamples$ is the set of samples assigned for collection from $\IndividualCentre$ by a selected route and $\ExpiredSamplesAt_\CurrentTime(\IndividualCentre) \subseteq \IndividualCentre_\CentreSamples$ is the set of samples at $\IndividualCentre$ whose deadlines are too close to allow collection and on-time delivery by any feasible primary-fleet route in a subsequent epoch. The derived attributes $\IndividualCentre_\CentreVolume^{\texttt{Post}}$ and $\IndividualCentre_\CentreDeadline^{\texttt{Post}}$ are recomputed from the updated sample set as summary features.

The exogenous step then incorporates the newly arrived samples $\ExogenousInformation_{\CurrentTime + 1}$ into their respective origin centres. Because vehicles are unaffected by sample arrivals, the fleet state carries over unchanged: $\IndividualVehicle_\VehicleReturn^{\CurrentTime + 1} = \IndividualVehicle_\VehicleReturn^{\texttt{Post}}$ for every $\IndividualVehicle \in \VehiclesSet$. For each centre $\IndividualCentre \in \CentresSet$, the pre-decision sample set at epoch $\CurrentTime + 1$ is obtained by adding each arriving sample whose origin is $\IndividualCentre$:
\begin{equation}
\IndividualCentre_\CentreSamples^{\CurrentTime + 1} = \IndividualCentre_\CentreSamples^{\texttt{Post}} \cup \{\IndividualSample \in \ExogenousInformation_{\CurrentTime + 1} : \IndividualSample_\SampleOrigin = \IndividualCentre\},
\label{eq:centre-exo-transition}
\end{equation}
with $\IndividualCentre_\CentreVolume^{\CurrentTime + 1}$ and $\IndividualCentre_\CentreDeadline^{\CurrentTime + 1}$ recomputed from $\IndividualCentre_\CentreSamples^{\CurrentTime + 1}$. The resulting pre-decision state $\StateNote_{\CurrentTime + 1} = (\VehicleState_{\CurrentTime + 1}, \CentreState_{\CurrentTime + 1})$ is then available to the dispatcher at the next decision epoch.

\subsection{Objective Function}

The dispatcher seeks a policy $\pi \in \Pi$, where $\Pi$ denotes the set of non-anticipative feasible policies mapping each state $\StateNote_\CurrentTime$ to a feasible decision vector $\DecisionsVector_\CurrentTime \in \FeasibleDecisions_\CurrentTime(\StateNote_\CurrentTime)$, that maximizes the expected total volume of blood samples served by the primary fleet over the operating horizon. Formally,
\begin{equation}
\max_{\pi \in \Pi} \; \mathbb{E}_{\ExogenousInformation = (\ExogenousInformation_1, \ldots, \ExogenousInformation_T)} \left[\; \sum_{\CurrentTime \in \Time} \Reward_\CurrentTime\!\left( \DecisionsVector_\CurrentTime^\pi\!\left(\StateNote_\CurrentTime^\pi(\ExogenousInformation)\right) \right) \;\middle|\; \StateNote_0 \right],
\label{eq:overall-objective}
\end{equation}
where $\DecisionsVector_\CurrentTime^\pi$ denotes the decision vector selected by policy $\pi$ in state $\StateNote_\CurrentTime^\pi(\ExogenousInformation)$, $\StateNote_0$ denotes the initial state in which all vehicles are idle at the hospital and all collection centres are empty, and the realized state trajectory evolves according to the recursion
\vspace{-1.2em}
\begin{subequations}
\label{eq:state-recursion}
\begin{align}
& \StateNote_0^\pi(\ExogenousInformation) = \StateNote_0, \\
& \StateNote_{\CurrentTime + 1}^\pi(\ExogenousInformation) = \SecondTransition\!\left( \FirstTransition\!\left(\StateNote_\CurrentTime^\pi, \DecisionsVector_\CurrentTime^\pi\!\left(\StateNote_\CurrentTime^\pi\right) \right), \ExogenousInformation_{\CurrentTime + 1} \right), \qquad \CurrentTime = 0, \ldots, T-1.
\end{align}
\end{subequations}
Because the random information vector $\ExogenousInformation$ is exogenous, the expectation in~\eqref{eq:overall-objective} is taken with respect to a distribution that does not depend on the policy itself.

The optimal value function satisfies the Bellman optimality equation:
\begin{equation}
V_\CurrentTime(\StateNote_\CurrentTime) = \max_{\DecisionsVector_\CurrentTime \in \FeasibleDecisions_\CurrentTime(\StateNote_\CurrentTime)} \left\{ \Reward_\CurrentTime(\DecisionsVector_\CurrentTime) + \DiscountRate \cdot \mathbb{E}_{\ExogenousInformation_{\CurrentTime + 1}} \!\left[ V_{\CurrentTime + 1}(\StateNote_{\CurrentTime + 1}) \right] \right\},
\label{eq:bellman}
\end{equation}
with terminal condition $V_T(\StateNote_T) = 0$, where $\StateNote_{\CurrentTime + 1} = \SecondTransition(\FirstTransition(\StateNote_\CurrentTime, \DecisionsVector_\CurrentTime), \ExogenousInformation_{\CurrentTime + 1})$ and $V_\CurrentTime(\StateNote_\CurrentTime)$ denotes the optimal expected cumulative reward obtainable from state $\StateNote_\CurrentTime$ onward. Since decisions are made over a finite single-day operating horizon, we adopt a discount factor of $\DiscountRate = 1$.

\section{Solution Methodology} \label{solutionmethodology}

This section develops a NeurADP approach tailored to the blood sample collection and delivery problem. We first motivate the approximation strategy by reformulating the Bellman optimality equation around the post-decision state. We then describe the two-step procedure that combines feasible-action enumeration with an integer linear matching problem, and finally present the dual vehicle-and-centre VFA and its training procedure.

\subsection{ADP and Post-Decision Value Approximation}

Although the Bellman optimality equation~\eqref{eq:bellman} admits an exact solution in principle via dynamic programming, direct solution is computationally intractable. At each epoch, candidate dispatch actions are generated by enumerating centre subsets and their route orderings, and these actions must be combined into a joint assignment across all vehicles and centres; the number of feasible decisions therefore grows exponentially with problem size. The state space is likewise prohibitively large, as it records each vehicle's remaining return time and the sample-level composition, including the volumes and deadlines, of the outstanding backlog at every centre. Moreover, current decisions determine which samples carry over to future epochs, while future arrivals are stochastic. Exact evaluation of the downstream value $\mathbb{E}[V_{\CurrentTime + 1}(\StateNote_{\CurrentTime + 1})]$ therefore requires accounting for all possible future arrival and dispatch trajectories.

To circumvent this intractability, we follow the ADP paradigm of \citet{powell2007approximate} and reformulate the Bellman optimality equation around the post-decision state $\PostDecisionState_\CurrentTime$ introduced in Section~\ref{modelformulation}. Recalling that the transition from $\StateNote_\CurrentTime$ to $\StateNote_{\CurrentTime + 1}$ decomposes into a deterministic post-decision step $\FirstTransition$ followed by a stochastic exogenous step $\SecondTransition$, the Bellman equation can equivalently be written in terms of a post-decision value function:
\label{eq:bellman-decomp}
\begin{align}
V_\CurrentTime(\StateNote_\CurrentTime) &= \max_{\DecisionsVector_\CurrentTime \in \FeasibleDecisions_\CurrentTime(\StateNote_\CurrentTime)} \Big\{ \Reward_\CurrentTime(\DecisionsVector_\CurrentTime) + V_\CurrentTime^{\texttt{Post}}(\PostDecisionState_\CurrentTime) \Big\}, \label{eq:bellman-pre} \\
V_\CurrentTime^{\texttt{Post}}(\PostDecisionState_\CurrentTime) &= \mathbb{E}_{\ExogenousInformation_{\CurrentTime + 1}}\!\left[\, V_{\CurrentTime + 1}(\StateNote_{\CurrentTime + 1}) \,\middle|\, \PostDecisionState_\CurrentTime \,\right]. \label{eq:bellman-post}
\end{align}
The key benefit of this reformulation is that the maximum over actions in~\eqref{eq:bellman-pre} is now over a deterministic quantity, because the expectation over exogenous arrivals has been absorbed into the post-decision value function~\eqref{eq:bellman-post}. Given a sufficiently accurate approximation of $V_\CurrentTime^{\texttt{Post}}$, the optimal action at each epoch could be obtained by solving a purely deterministic optimization problem rather than integrating over the exponentially large space of future arrival trajectories. NeurADP, introduced by \citet{shah2020neural}, takes exactly this approach: the post-decision value function is approximated by parameterized neural networks whose weights are learned by repeatedly simulating the system and backpropagating the temporal-difference error from each transition. We adapt this framework to the blood sample collection and delivery problem in the remainder of this section, with a dual value function decomposition tailored to the multi-entity structure of the problem as motivated in Section~\ref{sec:vfdecomp}.

\subsection{NeurADP for Blood Sample Collection and Delivery}

The NeurADP procedure at each decision epoch proceeds in two steps: (i) enumerate the feasible elementary actions at the current state, and (ii) solve an integer linear matching problem that selects the joint assignment maximizing the sum of immediate rewards plus decomposed post-decision value estimates. We employ two independently trained neural networks for these value estimates, one for vehicles and one for collection centres, whose architectures and feature inputs are detailed in Section~\ref{sec:nn_architecture}. 
We refer to the resulting model as \CombinedNeurADP.

\subsubsection{Two-Step Decomposition}
\label{sec:twostep}

At each decision epoch $\CurrentTime$, given the current state $\StateNote_\CurrentTime = (\VehicleState_\CurrentTime, \CentreState_\CurrentTime)$, we begin by enumerating the \revb{elementary action set $\FeasibleSet_\CurrentTime(\StateNote_\CurrentTime)$ defined in Section~\ref{modelformulation}}. For each candidate dispatch action, the optimal route $\OptimalRoute_{\IndividualVehicle, \CentreSubset}$ is determined by solving the shortest-path problem over permutations of $\CentreSubset$ per~\eqref{eq:tsp-route}, the greedy earliest-deadline-first sample collection yields $\UsedSamples_{\IndividualVehicle, \CentreSubset}$ and the associated collected volume $\CollectedVolume_{\IndividualVehicle, \CentreSubset}$, and the action is retained only if both the capacity constraint~\eqref{eq:capacity-constraint} and the deadline constraint~\eqref{eq:deadline-constraint} are satisfied. The feasible action set $\FeasibleSet_\CurrentTime(\StateNote_\CurrentTime)$ can be enumerated efficiently for the problem sizes we consider, in line with prior NeurADP applications to bounded-capacity dispatching \citep{shah2020neural,dehghan2025enhanced}.

Once the feasible action set is constructed, a binary decision variable \revb{$x_{\CurrentTime, (\IndividualVehicle, \CentreSubset)} \in \{0, 1\}$} is introduced for each action \revb{$(\IndividualVehicle, \CentreSubset) \in \FeasibleSet_\CurrentTime(\StateNote_\CurrentTime)$}, and the dispatcher solves the following integer linear matching problem:
\vspace{-1.9em}
\begin{subequations}
\label{eq:matching-ip}
\begin{align}
\texttt{MatchingIP}: \quad \max \ & \sum_{(\IndividualVehicle, \CentreSubset) \in \FeasibleSet_\CurrentTime(\StateNote_\CurrentTime)} x_{\CurrentTime, (\IndividualVehicle, \CentreSubset)} \cdot \ActionScore_\CurrentTime(\IndividualVehicle, \CentreSubset) \label{eq:match-obj} \\
\text{s.t.} \ & \eqref{eq:vehicle-constraint}, \eqref{eq:centre-constraint},
\end{align}
\end{subequations}
where \revb{$\ActionScore_\CurrentTime(\IndividualVehicle, \CentreSubset)$} denotes the total value of action \revb{$(\IndividualVehicle, \CentreSubset)$}, combining its immediate reward \revb{$\ImmediateReward(\IndividualVehicle, \CentreSubset)$} and its estimated downstream contribution obtained from the learned post-decision value functions described in the next subsection. Constraints~\eqref{eq:vehicle-constraint}--\eqref{eq:centre-constraint} enforce that each vehicle and each centre is covered by exactly one elementary action. The optimal solution of \texttt{MatchingIP} yields the decision vector $\DecisionsVector_\CurrentTime$ executed at the current epoch, after which the system advances according to the transition function of Section~\ref{modelformulation}.

\subsubsection{Dual Vehicle-and-Centre Value Function Approximation}
\label{sec:dualvfa}

To construct the scores \revb{$\ActionScore_\CurrentTime(\IndividualVehicle, \CentreSubset)$} required by \texttt{MatchingIP}, we approximate the post-decision value function $V_\CurrentTime^{\texttt{Post}}$ using a decomposition across entity types, namely as a sum of per-entity contributions, with one term per vehicle and one term per collection centre:
\begin{equation}
\label{eq:dual-decomposition}
V_\CurrentTime^{\texttt{Post}}(\PostDecisionState_\CurrentTime) \approx \sum_{\IndividualVehicle \in \VehiclesSet} \PostValueVeh\!\left(\IndividualVehicle^{\texttt{Post}}, \AuxFeatures_\CurrentTime\right) + \sum_{\IndividualCentre \in \CentresSet} \PostValueCtr\!\left(\IndividualCentre^{\texttt{Post}}, \AuxFeatures_\CurrentTime\right),
\end{equation}
where $\PostValueVeh$ and $\PostValueCtr$ are two parameterized neural networks producing scalar value estimates, $\IndividualVehicle^{\texttt{Post}}$ and $\IndividualCentre^{\texttt{Post}}$ are the post-decision attribute tuples of individual vehicles and centres respectively, and $\AuxFeatures_\CurrentTime$ is a vector of global auxiliary features summarizing the system-wide state at epoch $\CurrentTime$, shared across both networks as context. Each network thus conditions on its own entity's post-decision state and on a compact summary of the surrounding system, enabling it to learn entity-specific value dynamics while remaining aware of the global operational context. The specific architecture, input feature set, and hyperparameter choices are reported in Section~\ref{sec:nn_architecture}.

Given this decomposition, \revb{the score of a dispatch action $(\IndividualVehicle, \CentreSubset) \in \FeasibleDispatch_\CurrentTime(\StateNote_\CurrentTime)$} is
\begin{equation}
\label{eq:action-score}
\ActionScore_\CurrentTime(\IndividualVehicle, \CentreSubset) = \ImmediateReward(\IndividualVehicle, \CentreSubset) + \DiscountRate \cdot \Bigg[\, \PostValueVeh\!\left(\IndividualVehicle^{\texttt{Post}}, \AuxFeatures_\CurrentTime\right) + \sum_{\IndividualCentre \in \CentreSubset} \PostValueCtr\!\left(\IndividualCentre^{\texttt{Post}}, \AuxFeatures_\CurrentTime\right) \,\Bigg],
\end{equation}
where $\IndividualVehicle^{\texttt{Post}}$ and $\IndividualCentre^{\texttt{Post}}$ are computed by applying the action deterministically to $\StateNote_\CurrentTime$ via the post-decision transition equations~\eqref{eq:vehicle-transition} and~\eqref{eq:centre-post-transition}. \revb{The null actions carry no immediate reward and are scored by their single post-decision term, $\ActionScore_\CurrentTime(\IndividualVehicle, \NullAction) = \DiscountRate \, \PostValueVeh(\IndividualVehicle^{\texttt{Post}}, \AuxFeatures_\CurrentTime)$ and $\ActionScore_\CurrentTime(\NoVehicle, \{\IndividualCentre\}) = \DiscountRate \, \PostValueCtr(\IndividualCentre^{\texttt{Post}}, \AuxFeatures_\CurrentTime)$.} Because each centre and each vehicle is evaluated independently by its own network, the cost of scoring all feasible actions scales linearly with the number of entities involved rather than exponentially with the size of the joint state, making the score computation tractable even as the system grows.

Both value networks are trained jointly via prioritized experience replay and soft target updates, following standard deep reinforcement learning practice. During training, the system simulates a day of operation: at every decision point, the feasible action set is constructed, actions are scored using the current \emph{online} networks, and the resulting matching problem is solved after adding a decaying Gaussian perturbation to each action's score, so that the policy explores diverse dispatch patterns early in training and gradually shifts toward exploitation as training progresses. The realized transition, together with the pre-decision state, the post-decision states of all entities, and the feasible action set, is stored as an experience in a prioritized replay buffer in which experiences with higher temporal-difference error are sampled more frequently. After a fixed number of decision epochs, a minibatch of experiences is drawn from the buffer, and each network is updated by minimizing the mean squared error between its predicted post-decision values and supervised targets. The supervised targets themselves are constructed by re-solving \texttt{MatchingIP} on the stored feasible action set using slowly tracking \emph{target networks}, whose weights are updated after each minibatch via a soft Polyak interpolation toward the online network weights. This target network mechanism stabilizes the bootstrapped updates that are characteristic of temporal-difference learning. The detailed neural network architectures and the full training procedure (Algorithm~\ref{alg:training}) are provided in Appendix~\ref{ec:nnarch}.

\section{Experimental Setup} \label{experimentalsetup}
All methods are implemented in Python 3.10.9 and run on a Google Cloud e2-standard-8 virtual machine equipped with 8 vCPUs, 32 GB of memory, and an AMD Rome CPU (x86/64). Integer linear programming (ILP) models are solved using IBM ILOG CPLEX Optimization Studio version 22.1.1. Below, we present the experimental setting used to evaluate the proposed policies for the blood sample collection and delivery problem, including the collection network, stochastic sample-arrival process, neural network architectures, and benchmark policies.

\subsection{Dataset Description}\label{sec:dataset}

We evaluate the policies on a collection network based on the Greater Toronto Area, comprising 32 licensed specimen collection centres and a single central hospital (Toronto General Hospital) that serves as both the depot and the central laboratory. Centre locations are drawn from real downtown-Toronto collection facilities, and each centre is assigned a relative demand prevalence---grouped into high-, medium-, and low-volume tiers---derived from publicly available review volume used as a proxy for patient traffic. Travel times between locations use Haversine distances at an effective courier speed of 20~km/h, plus a 5-minute per-centre service time for sample handoff and loading. Blood samples arrive according to a non-homogeneous Poisson process calibrated to empirical outpatient clinic data \citep{feng2023dataset}, producing a bimodal arrival pattern with morning and afternoon peaks; arrivals are distributed across centres with stochastic spatial variation, and each sample is assigned an integer volume drawn uniformly from $\{1,2,3,4,5\}$ units. The network topology and the relative demand intensities are grounded in publicly available Toronto data, whereas the travel-time model, the review-volume proxy, and the spatial and volume distributions are modeling approximations rather than measurements from an operational system. The full network construction, the demand-calibration details, and the corresponding network map and arrival-rate figures are provided in Appendix~\ref{ec:datasetdetail}.

\subsection{Neural Network Architectures and Training}\label{sec:nn_architecture}

As described in Section~\ref{solutionmethodology}, the \CombinedNeurADP policy approximates the dual post-decision value function with two feedforward neural networks---one for vehicles and one for collection centres---that share a common backbone of fully connected layers (three layers of 300 hidden units) with Exponential Linear Unit (ELU) activations and differ in their entity-specific inputs. The vehicle network conditions on the vehicle's normalized time until return; the centre network conditions on a learnable location embedding and the centre's stored volume, earliest deadline, and sample count. Both networks additionally receive five system-level auxiliary features, namely, the number of incoming samples, the number of available vehicles, the average deadline urgency across stored samples, the average stored volume per centre, and the average vehicle return time, together with a learnable embedding of the current epoch, which supply the shared global context.
The complete network architectures, input-feature definitions, training hyperparameters, and the training algorithm are provided in Appendix~\ref{ec:nnarch}.

\subsection{Benchmark Policies}\label{sec:benchmarks}

We compare \CombinedNeurADP with four benchmark policies that isolate the role of each value-function component and the value of non-myopic dispatching.
\begin{itemize}[leftmargin=1.2em]\setItemSep{-0.2em}
\item \VehicleNeurADP: This ablation uses the same NeurADP training and dispatching framework as \CombinedNeurADP but retains only the vehicle value function, i.e., only the vehicle network is trained and dispatch is guided solely by vehicle-level post-decision values, with the centre component set to zero, isolating the contribution of vehicle-level foresight.

\item \CentreNeurADP: This ablation is the symmetric counterpart, retaining only the centre value function, i.e., only the centre network is trained and dispatch is guided solely by centre-level post-decision values. Together with \VehicleNeurADP, it separates centre-level awareness from vehicle-level foresight.

\item \ILPThroughput: This myopic policy solves an ILP at each decision epoch to maximize the ratio of collected volume to travel time, without incorporating any lookahead. By favouring routes with high collection efficiency and rapid vehicle turnover, it represents a greedy strategy focused on throughput.

\item \ILPVolume: This myopic policy also solves an ILP at each epoch, but instead maximizes the total volume collected per dispatch. It favours routes that consolidate pickups to utilize vehicle capacity more fully, representing a greedy strategy that prioritizes immediate volume over routing efficiency.

\end{itemize}


\subsection{Baseline Configuration and Evaluation Protocol}\label{sec:baseline}
Our experimental setting follows the two-tier operational model of \citet{zufferey2016dynamic}, developed with a Geneva blood sample collection laboratory. There, a primary fleet of scooters performs most collections while samples that cannot be served on time are diverted to external couriers or taxis at higher cost; scooters travel about 10\% faster than cars in dense urban traffic, suiting them to the strict 90-minute viability windows, and the laboratory tolerates up to 10\% of requests being handled by the secondary channel. We adopt an analogous context for the Toronto network: a primary fleet of mopeds or bikes performs real-time collections, and any unserved samples are fulfilled by external vehicles at higher cost. Accordingly, we optimize the primary fleet's dispatch decisions to maximize the share of sample volume delivered within the deadline (i.e., the \emph{service rate}), which directly reduces dependence on the costly secondary channel; a service rate above 90\% is a strong operational target, consistent with the levels reported in practice.

Unless stated otherwise, we use the baseline configuration in Table~\ref{table:baseline}. Time is discretized into $\EpochLength = 5$-minute decision epochs; the active arrival period spans 8.5 hours and is followed by a terminal clearance period of $\DeadPeriod = 90$ minutes during which no new samples arrive but vehicles continue to collect the outstanding backlog. Under this configuration the system receives, on average, 356.7 samples per day, or approximately 1{,}067 volume units. Each simulated day is a single demand sample path; every policy is evaluated on the same 20 independent test days. We study sensitivity by varying one parameter at a time, namely, sample deadline ($\Deadline$), vehicle capacity ($\Capacity$), fleet size ($\NumAgents$), and maximum locations per trip ($\LocsToVisit$), while holding the others at their baseline values.

\begin{table}[!ht]
\centering
\caption{Baseline configuration, used in all experiments unless otherwise stated.}
\label{table:baseline}
\setlength{\tabcolsep}{6pt}
\renewcommand{\arraystretch}{0.91}
\resizebox{0.5\textwidth}{!}{
\begin{tabular}{llc}
\toprule
\textbf{Parameter} & \textbf{Symbol} & \textbf{Value} \\
\midrule
Fleet size & $\NumAgents$ & 4 vehicles \\
Vehicle capacity & $\Capacity$ & 15 units \\
Sample deadline & $\Deadline$ & 90 minutes \\
Maximum centres per trip & $\LocsToVisit$ & 2 \\
Decision-epoch length & $\EpochLength$ & 5 minutes \\
Terminal clearance period & $\DeadPeriod$ & 90 minutes \\
Active arrival period & --- & 8.5 hours \\
Test days evaluated & --- & 20 \\
\bottomrule
\end{tabular}
}
\end{table}

\section{Numerical Results}\label{results}

This section evaluates \CombinedNeurADP against all benchmark policies across several operational settings.
Each sensitivity table reports the following metrics. The ``Combined'' column gives the service rate, measured as the percentage of sample volume served by \CombinedNeurADP and reported as mean $\pm$ standard error over the 20 test days. The remaining columns, ``\% Over Vehicle'', ``\% Over Centre'', ``\% Over ILP-Thr.'', and ``\% Over ILP-Vol.'', report the percentage improvement of \CombinedNeurADP relative to each respective benchmark, computed as the absolute difference in mean service rate. A positive value indicates that \CombinedNeurADP exceeds the benchmark by that margin.


We analyze the sample-deadline experiment in full in the main text as a representative case, reporting its headline service rates together with the underlying operational trade-offs. For the remaining three experiments we report the headline service rates here and defer the detailed operational analysis---dispatch frequency, route composition, average travel time, and spatial-coverage figures---to Appendix~\ref{ec:opanalysis}.

\subsection{Impact of Sample Deadline}
We first examine how varying the sample deadline ($\Deadline$) affects system performance (Table~\ref{table:deadline_table}). Relaxing the deadline improves outcomes across all policies, with service rates increasing from roughly 82--87\% at 60 minutes to 94--97\% at 120 minutes. \CombinedNeurADP achieves the highest service rate at all three deadline levels, although its advantage over the baselines narrows as deadlines become less restrictive, reflecting the reduced difficulty of the dispatch problem.

\begin{table}[!ht]
\centering
\caption{Impact of sample deadline on volume served (\%). Results averaged over 20 test days. ``\% Over'' columns show percentage-point improvement of \CombinedNeurADP over each policy.}
\label{table:deadline_table}
\setlength{\tabcolsep}{5pt}
\renewcommand{\arraystretch}{1.01}
\resizebox{0.9\textwidth}{!}{
\begin{tabular}{lccccc}
\toprule
\textbf{Deadline} & \textbf{Combined} & \textbf{\%  Over Vehicle } & \textbf{\%  Over Centre } & \textbf{\%  Over ILP-Thr. } & \textbf{\%  Over ILP-Vol. } \\
\midrule
60 min & $86.83 \pm 0.43$ & +2.24 & +1.36 & +3.11 & +4.59 \\
90 min & $92.96 \pm 0.24$ & +1.04 & +1.35 & +2.23 & +3.31 \\
120 min & $96.73 \pm 0.21$ & +0.81 & +1.78 & +1.53 & +2.59 \\
\bottomrule
\end{tabular}
}
\end{table}

\begin{figure}[!ht]
\centering
\includegraphics[width=0.49\textwidth]{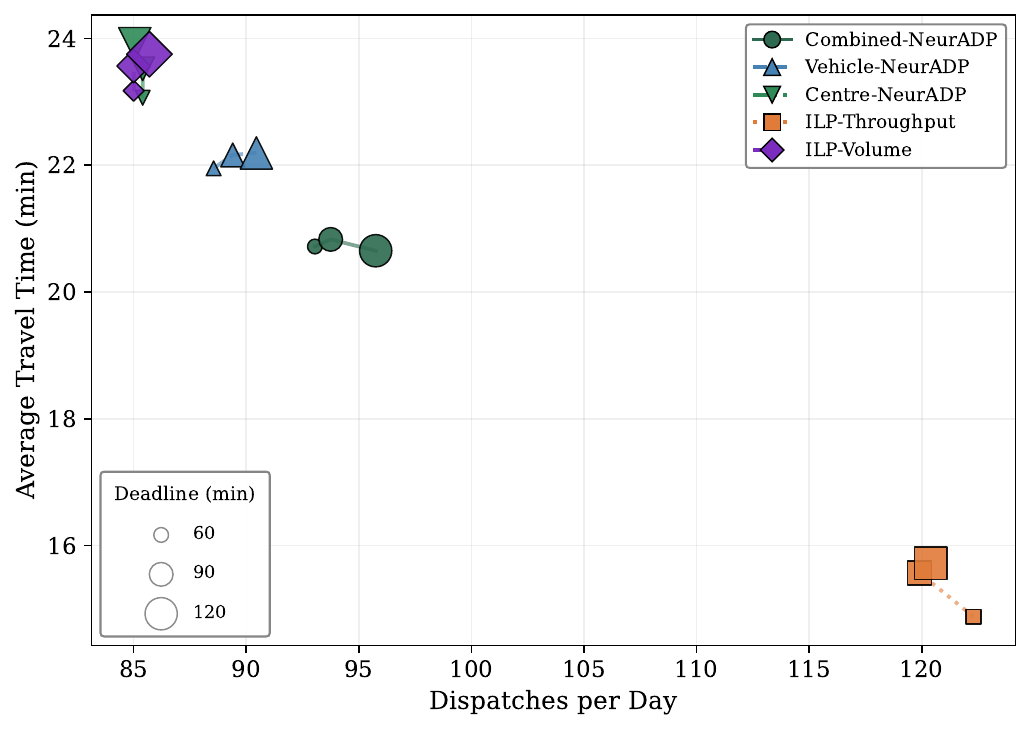}
\caption{Dispatch frequency vs.\ average travel time across policies under varying sample deadlines. Marker size indicates the deadline.}
\label{fig:deadline_tradeoff}
\end{figure}

Figure~\ref{fig:deadline_tradeoff} plots dispatch frequency against average travel time across policies. We observe a clear structural trade-off: policies that dispatch more frequently tend to construct shorter routes, whereas policies that consolidate more aggressively incur longer travel times. \CombinedNeurADP occupies a balanced region of this trade-off, maintaining moderate dispatch frequency and relatively short travel times. Its position remains stable across the deadline settings, indicating that it does not rely on reactive changes in route length. Instead, it adapts primarily through centre selection, which allows it to maintain both timely collection and spatial coverage as constraints tighten or relax.
The myopic baselines lie at opposite extremes. \ILPThroughput operates in a high-dispatch, low-travel-time regime, fragmenting routes into short, single-centre trips; this improves turnaround but leads to poor spatial coverage, with some centres receiving disproportionately few visits and experiencing frequent expirations. \ILPVolume adopts a low-dispatch, high-travel-time strategy that consistently favours multi-centre consolidation; this balances coverage but delays collection, causing samples to expire while waiting to be bundled into routes. The ablation variants fall between these extremes but do not reach the same balance as \CombinedNeurADP. \VehicleNeurADP operates closer to the consolidated regime, dispatching less frequently and incurring longer travel times, whereas \CentreNeurADP pushes consolidation further and closely resembles \ILPVolume. Taken together, these results show that deadline performance is driven by how effectively a policy navigates the trade-off between dispatch frequency and route duration: by integrating vehicle- and centre-level value functions, \CombinedNeurADP sustains frequent, efficient dispatches while preserving balanced spatial coverage.

\subsection{Impact of Vehicle Capacity}
We next examine how vehicle capacity ($\Capacity$) affects system performance (Table~\ref{table:capacity_table}). Increasing capacity from 10 to 15 yields noticeable gains across all policies, with \CombinedNeurADP improving from 87.84\% to 92.96\%. Further increasing capacity to 20 provides only marginal improvement for \CombinedNeurADP (93.71\%), indicating diminishing returns in the tested setting. This pattern suggests that, under baseline demand, capacity is no longer the main binding constraint once it exceeds approximately 15 units.

\begin{table}[!ht]
\centering
\caption{Impact of vehicle capacity on volume served (\%). Results averaged over 20 test days. ``\% Over'' columns show percentage-point improvement of \CombinedNeurADP over each policy.}
\label{table:capacity_table}
\setlength{\tabcolsep}{5pt}
\renewcommand{\arraystretch}{1.01}
\resizebox{0.9\textwidth}{!}{
\begin{tabular}{lccccc}
\toprule
\textbf{Capacity} & \textbf{Combined} & \textbf{\%  Over Vehicle } & \textbf{\%  Over Centre } & \textbf{\%  Over ILP-Thr. } & \textbf{\%  Over ILP-Vol. } \\
\midrule
10 & $87.84 \pm 0.59$ & +2.73 & +9.41 & +2.65 & +8.62 \\
15 & $92.96 \pm 0.24$ & +1.04 & +1.35 & +2.23 & +3.31 \\
20 & $93.71 \pm 0.28$ & +1.52 & +2.10 & +2.10 & +3.14 \\
\bottomrule
\end{tabular}
}
\end{table}

\subsection{Impact of Fleet Size}

We next examine how fleet size ($\NumAgents$) affects system performance (Table~\ref{table:agents_table}). This parameter has the largest impact of any considered: \CombinedNeurADP improves from 84.49\% with 3 vehicles to 96.81\% with 5 vehicles, a gain of more than 12 percentage points. As in previous experiments, the performance gap between policies is widest under scarcity and narrows as resources increase.

\begin{table}[!ht]
\centering
\caption{Impact of fleet size on volume served (\%). Results averaged over 20 test days. ``\% Over'' columns show percentage-point improvement of \CombinedNeurADP over each policy.}
\label{table:agents_table}
\setlength{\tabcolsep}{5pt}
\renewcommand{\arraystretch}{1.01}
\resizebox{0.9\textwidth}{!}{
\begin{tabular}{lccccc}
\toprule
\textbf{Agents} & \textbf{Combined} & \textbf{\%  Over Vehicle } & \textbf{\%  Over Centre } & \textbf{\%  Over ILP-Thr. } & \textbf{\%  Over ILP-Vol. } \\
\midrule
3 & $84.49 \pm 0.47$ & +2.28 & +4.52 & +2.79 & +6.94 \\
4 & $92.96 \pm 0.24$ & +1.04 & +1.35 & +2.23 & +3.31 \\
5 & $96.81 \pm 0.18$ & +0.69 & +0.65 & +1.09 & +1.56 \\
\bottomrule
\end{tabular}
}
\end{table}

\subsection{Impact of Maximum Locations Per Trip}
Finally, we examine the effect of varying the maximum number of centres a vehicle may visit per trip ($\LocsToVisit$). This parameter directly controls the complexity of the routing decision: at $\LocsToVisit = 1$, the problem reduces to selecting a single centre, whereas larger values require jointly selecting and sequencing multiple stops. Table~\ref{table:locations_table} presents the results.

\begin{table}[!ht]
\centering
\caption{Impact of locations per trip on volume served (\%). Results averaged over 20 test days. ``\% Over'' columns show percentage-point improvement of \CombinedNeurADP over each policy; a negative value indicates that the benchmark outperforms \CombinedNeurADP.}
\label{table:locations_table}
\setlength{\tabcolsep}{5pt}
\renewcommand{\arraystretch}{1.01}
\resizebox{0.9\textwidth}{!}{
\begin{tabular}{lccccc}
\toprule
\textbf{Locations} & \textbf{Combined} & \textbf{\%  Over Vehicle } & \textbf{\%  Over Centre } & \textbf{\%  Over ILP-Thr. } & \textbf{\%  Over ILP-Vol. } \\
\midrule
1 & $88.76 \pm 0.24$ & +1.30 & $-$0.66 & +2.67 & +1.84 \\
2 & $92.96 \pm 0.24$ & +1.04 & +1.35 & +2.23 & +3.31 \\
3 & $93.58 \pm 0.36$ & +2.05 & +8.63 & +1.96 & +8.64 \\
\bottomrule
\end{tabular}
}
\end{table}

The relative ordering of policies aligns with the other experiments: \CombinedNeurADP attains the highest service rate at every setting except $\LocsToVisit = 1$, where routing is trivial and \CentreNeurADP edges it out by $0.66$ percentage points because centre selection alone then suffices. As the action space widens to $\LocsToVisit = 3$, the gap over the ablation variants grows, since the vehicle component is needed to curb excessive route length while the centre component preserves spatial coverage. The detailed route-composition analysis is provided in Appendix~\ref{ec:opanalysis}.

\section{Conclusion}\label{conclusion}

This paper develops a NeurADP framework for real-time blood sample collection under stochastic demand and perishability constraints. The approach combines vehicle- and collection-centre value functions within an ILP dispatch model, allowing dispatch decisions to account for both current collection opportunities and future system impact. On a realistic Toronto-based network, \CombinedNeurADP outperforms myopic ILP baselines and single-component ablations, with the largest gains under tight deadlines, limited capacity, and small fleet sizes.

The main operational insight is that effective dispatching requires balancing route consolidation with timely collection. The dual value-function structure supports this balance by combining supply-side information about fleet availability with demand-side information about centre-level urgency. For healthcare logistics managers, the results indicate that foresighted dispatching is most valuable when resources are constrained, because better dispatch decisions can increase service levels and reduce reliance on costly backup options.
Future work could extend the framework to multi-depot networks, time-dependent travel times, and richer value-function mixing architectures, as well as related healthcare logistics applications such as pharmaceutical distribution and organ transport.

\bibliographystyle{apalike}
\bibliography{references}

\begin{APPENDICES}
\section{Collection Network Construction and Demand Calibration}\label{ec:datasetdetail}

This section provides the full construction of the semi-synthetic collection network and the demand-calibration details summarized in Section~\ref{sec:dataset} of the main paper.

We construct a Greater Toronto Area collection network consisting of 32 licensed specimen collection centres and a single central hospital. The collection centres are based on real-world specimen-collection facilities across downtown Toronto, including providers such as LifeLabs and Dynacare. The central hospital is Toronto General Hospital, one of the largest academic health science centres in Canada and a major hub for laboratory diagnostics, making it a natural choice as the central processing facility in our network. Figure~\ref{fig:location_map} illustrates the geographic distribution of the network, with marker sizes proportional to the arrival prevalence at each centre. To approximate realistic demand heterogeneity, each centre is assigned a relative prevalence based on publicly available review volume, used as a proxy for patient traffic. Centres are then grouped into high-, medium-, and low-volume tiers. Travel times between all location pairs are computed using Haversine distances and an average travel speed of 20 km/h, reflecting typical effective speeds for bicycle and moped couriers operating in dense urban environments with congestion and frequent stops. An additional 5-minute service time is included at each centre to account for sample handoff and loading. We adjust the 10- to 15-minute stopping time reported by \citet{grasas2014improvement} to reflect a more streamlined urban pickup process, where smaller vehicles reduce parking delays and lower-volume, more frequent pickups reduce handling time.

\begin{figure}[!ht]
\centering
\includegraphics[width=0.65\textwidth]{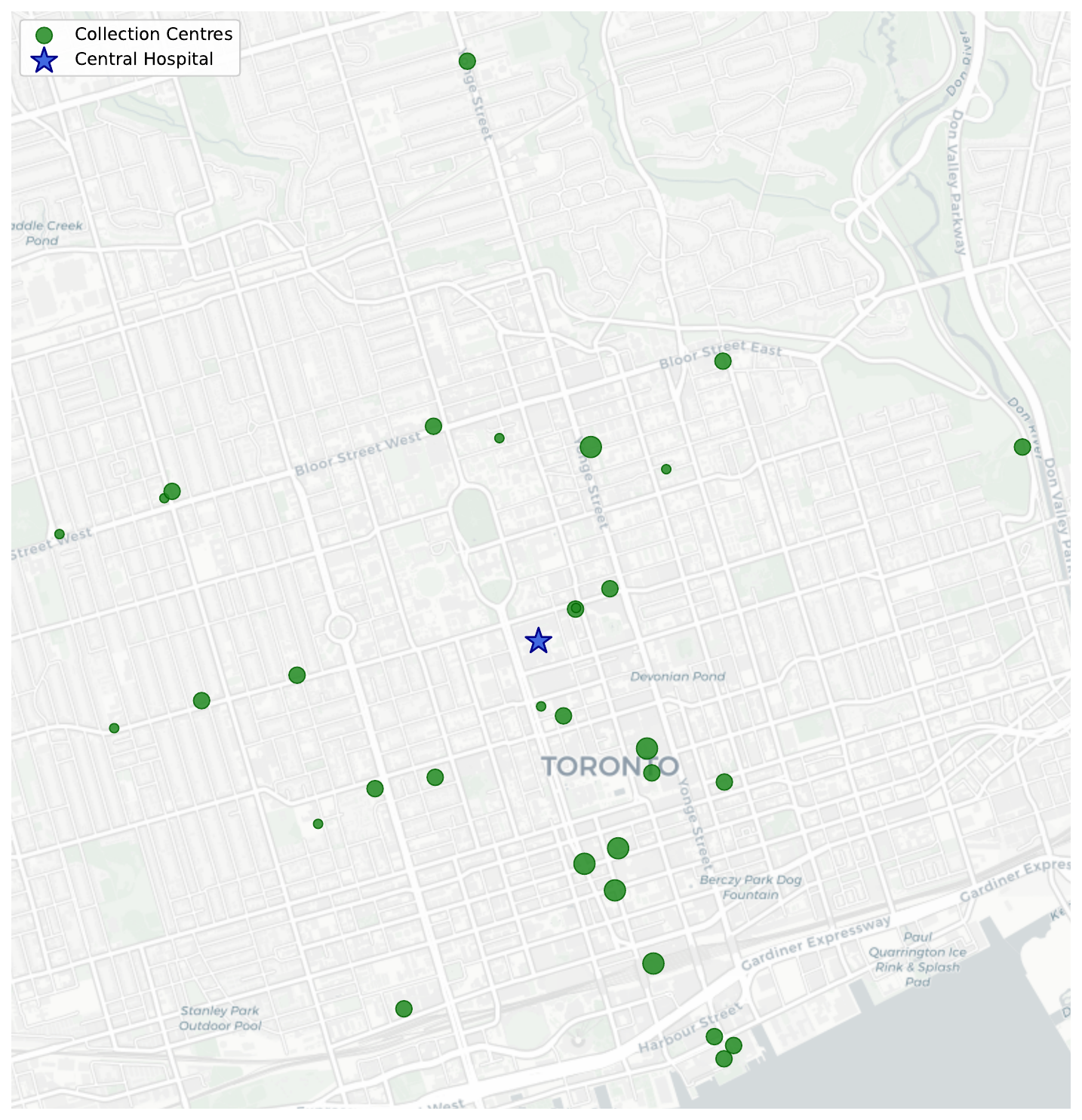}
\caption{Geographic distribution of collection centres (green) and hospital (blue star) in the Toronto network. Marker sizes reflect relative arrival prevalence.}
\label{fig:location_map}
\end{figure}

Blood samples arrive stochastically throughout the day according to a non-homogeneous Poisson process calibrated to empirical outpatient clinic data from \citet{feng2023dataset}. The resulting arrival pattern, shown in Figure~\ref{fig:arrival_dist}, exhibits a bimodal structure with morning and afternoon peaks separated by a midday dip, and is mapped to 5-minute epoch rates. To introduce spatial variation, the distribution of arrivals across centres at each epoch is sampled from a symmetric Dirichlet distribution with concentration parameter $\alpha = 2.0$, so that demand shifts stochastically across the network over the course of the day. Finally, each sample is assigned an integer volume drawn uniformly from $\{1, 2, 3, 4, 5\}$ units.

Taken together, these choices define a realistic testbed: the network topology and the relative demand intensities across centres are grounded in publicly available Toronto data, whereas the travel-time model (Haversine distances at a constant 20~km/h), the review-volume proxy for patient traffic, and the spatial and volume distributions are modeling approximations rather than measurements from an operational system.

\begin{figure}[!ht]
\centering
\includegraphics[width=0.805\textwidth]{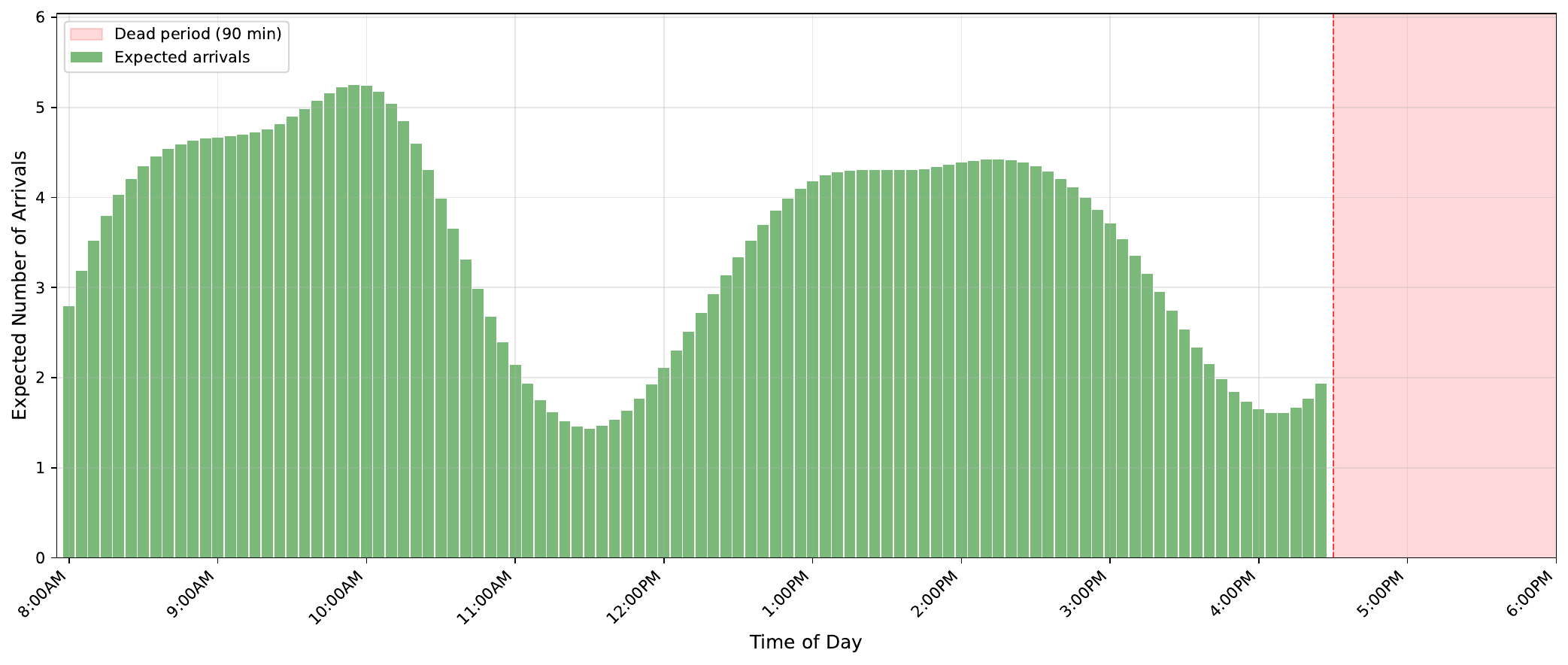}
\caption{Expected number of sample arrivals per 5-minute epoch throughout the operating day. The shaded red region indicates the 90-minute dead period during which no new samples arrive.}
\label{fig:arrival_dist}
\end{figure}

\section{Neural Network Architectures and Training Details}\label{ec:nnarch}

As described in Section~\ref{solutionmethodology} of the main paper, the \CombinedNeurADP policy decomposes the system-level value function into two components: one for vehicles and one for collection centres. Each is approximated by a separate feedforward neural network that maps a post-decision state to a scalar value estimate. Both networks share a common architectural backbone composed of fully connected layers with Exponential Linear Unit (ELU) activations, differing primarily in their input structure and the entity they represent.

The vehicle network estimates the downstream value of a vehicle's post-decision state. Its input consists of the vehicle's normalized time until return, along with five system-level auxiliary features that capture the current operational context: the number of incoming samples, the number of available vehicles, the average deadline urgency across stored samples, the average volume stored across centres, and the average vehicle return time. The current decision epoch is passed through a learnable time embedding to capture temporal patterns in the problem dynamics. This embedding is concatenated with the scalar features and passed through three fully connected layers of 300 hidden units each, producing a scalar value estimate. The centre network follows a similar architecture but additionally conditions on centre-specific state information. Each input includes a learnable location embedding that encodes the centre's identity, along with the centre's current stored volume, earliest sample deadline, and number of stored samples. These centre-specific features are concatenated with the same time embedding and system-level auxiliary features used by the vehicle network before being passed through the same three-layer feedforward architecture.

At each decision epoch, the total value of a candidate action is computed as the sum of its immediate volume reward and the downstream values of all entities involved: the vehicle value from the vehicle network and the individual centre values from the centre network for each centre visited. These total action values are then passed to the ILP, which selects the jointly optimal assignment of vehicles to centres by maximizing the sum of action values subject to the constraint that each vehicle and each centre is assigned to exactly one action.

Both networks are trained concurrently within a reinforcement learning framework using prioritized experience replay and soft target updates. During each decision epoch, the system interacts with a simulated environment to generate state observations, compute action values via the ILP, and execute dispatch decisions. These interactions are stored in a prioritized replay buffer, where experiences with higher temporal-difference error are sampled more frequently during training. Periodically, a batch of experiences is drawn and used to update both networks by minimizing the mean squared error between their predicted values and supervised targets. For each sampled experience, per-entity supervised targets are constructed by re-solving \texttt{MatchingIP} at the stored next state using the target networks, which yields a bootstrap assignment together with the resulting post-decision states of every vehicle and centre. The realized collection is allocated across the two value functions according to where each unit of volume is picked up. Writing $q_\CurrentTime(\IndividualCentre) = \sum_{\IndividualSample \in \UsedSamplesAt_\CurrentTime(\IndividualCentre)} \IndividualSample_\SampleVolume$ for the volume collected from centre $\IndividualCentre$ during the epoch, a vehicle dispatched to centre subset $\CentreSubset$ collects the route total $\CollectedVolume_{\IndividualVehicle,\CentreSubset} = \sum_{\IndividualCentre \in \CentreSubset} q_\CurrentTime(\IndividualCentre)$. The vehicle target is the total volume collected along the vehicle's assigned route plus the discounted target-network value of the vehicle's resulting post-decision state, and each centre target is that centre's own collected volume plus the discounted target-network value of the centre's resulting post-decision state:
\begin{equation}
y^{\texttt{Veh}}_{\IndividualVehicle} = \CollectedVolume_{\IndividualVehicle,\CentreSubset} + \DiscountRate\, \TargetPostValueVeh\!\left(\IndividualVehicle^{\texttt{Post}}\right),
\qquad
y^{\texttt{Ctr}}_{\IndividualCentre} = q_\CurrentTime(\IndividualCentre) + \DiscountRate\, \TargetPostValueCtr\!\left(\IndividualCentre^{\texttt{Post}}\right),
\label{ec:targets}
\end{equation}
where $\IndividualVehicle^{\texttt{Post}}$ and $\IndividualCentre^{\texttt{Post}}$ denote the post-decision states induced by the re-solved bootstrap assignment, and an idle vehicle or an unvisited centre takes a zero immediate term and is bootstrapped from its carried-over post-decision state. Because every collected unit of volume originates at exactly one centre and is transported by exactly one vehicle, this attribution is consistent with the additive decomposition in~\eqref{eq:dual-decomposition}: the vehicle head is trained on route-level collection, while the centre heads are trained on the matching per-centre collection. Separate slowly tracking target networks are maintained for both the vehicle and centre models to enhance training stability, and an adaptive learning rate schedule is employed for improved convergence behaviour. To encourage exploration during training, Gaussian noise with a decaying magnitude is added to the action values prior to the ILP optimization step, allowing the policy to discover high-quality dispatch strategies before gradually shifting toward exploitation. Together, these mechanisms enable both neural networks to learn temporally extended value functions under operational uncertainty, supporting high-quality vehicle dispatch decisions throughout the operating day. Consistent with the MDP timing convention of Section~\ref{modelformulation}, each epoch's dispatch decision is made on the pre-decision state $\StateNote_\CurrentTime$ before that epoch's arrivals are realized; the sampled arrivals enter only through the transition to the next pre-decision state, so the learned policy remains non-anticipative. The full training procedure is summarized in Algorithm~\ref{alg:training}, and the hyperparameter settings are listed in Table~\ref{table:hyperparams}.

\begin{algorithm}[!ht]
\caption{Training procedure for \CombinedNeurADP.}
\label{alg:training}
\begingroup
\setlength{\baselineskip}{0.9\baselineskip}
\begin{algorithmic}[1]
\State \textbf{Input:} number of training days $\Days$, minibatch size $|\mathcal{M}_{\text{batch}}|$, update frequency $K$, soft target-update rate $\StepSize$, exploration schedule
\State \textbf{Output:} trained vehicle and centre value networks $\PostValueVeh$ and $\PostValueCtr$
\State Initialize online networks $\PostValueVeh$, $\PostValueCtr$ and target networks $\TargetPostValueVeh \leftarrow \PostValueVeh$, $\TargetPostValueCtr \leftarrow \PostValueCtr$
\State Initialize empty prioritized replay buffer $\mathcal{B}$
\For{day $\SingleDay = 1, \ldots, \Days$}
    \State Initialize system state $\StateNote_0$ with all vehicles idle and all centres empty
    \For{epoch $\CurrentTime = 0, \ldots, T-1$}
        \State Compute global auxiliary features $\AuxFeatures_\CurrentTime$ from the current pre-decision state $\StateNote_\CurrentTime$
        \State Enumerate feasible actions $\FeasibleSet_\CurrentTime(\StateNote_\CurrentTime)$ and their post-decision states
        \State Score each action via~\eqref{eq:action-score} using the \emph{online} networks, adding a decaying Gaussian exploration term
        \State Solve \texttt{MatchingIP}~\eqref{eq:matching-ip} to obtain decision vector $\DecisionsVector_\CurrentTime$
        \State Execute $\DecisionsVector_\CurrentTime$, observe reward $\Reward_\CurrentTime(\DecisionsVector_\CurrentTime)$, and form the post-decision state $\PostDecisionState_\CurrentTime$
        \State Sample arrivals $\ExogenousInformation_{\CurrentTime+1}$ \emph{after} the decision, and obtain the next pre-decision state $\StateNote_{\CurrentTime+1}$
        \State Store experience $(\StateNote_\CurrentTime, \FeasibleSet_\CurrentTime, \PostDecisionState_\CurrentTime, \PostDecisionState_{\CurrentTime-1})$ in $\mathcal{B}$
        \If{$\CurrentTime \bmod K = 0$ \textbf{and} $|\mathcal{B}|$ exceeds warm-up threshold}
            \State Sample a minibatch of experiences from $\mathcal{B}$ with priority-weighted sampling
            \For{each experience in the minibatch}
                \State Re-score all stored feasible actions via~\eqref{eq:action-score} using the \emph{target} networks
                \State Re-solve \texttt{MatchingIP} to extract per-vehicle and per-centre supervised value targets
            \EndFor
            \State Update $\PostValueVeh$ and $\PostValueCtr$ by minimizing the mean squared error between predicted values and supervised targets, weighted by importance sampling corrections
            \State Update replay buffer priorities for sampled experiences using the new squared errors
            \State Soft-update target networks: $\theta^{\texttt{target}} \leftarrow \StepSize \cdot \theta + (1 - \StepSize) \cdot \theta^{\texttt{target}}$
        \EndIf
    \EndFor
\EndFor
\State \Return $\PostValueVeh$, $\PostValueCtr$
\end{algorithmic}
\endgroup
\end{algorithm}

\begin{table}[!ht]
\centering
\begingroup
\setlength{\tabcolsep}{2pt}
\renewcommand{\arraystretch}{0.93}

\caption{Hyperparameter settings for the \CombinedNeurADP{} training procedure
(Algorithm~\ref{alg:training}). Architectural, horizon, and discount settings
are used in all experiments; the remaining operational parameters follow the baseline configuration in Table~\ref{table:baseline} of the main paper}.
\label{table:hyperparams}

\begin{tabular}{@{}ll@{}}
\toprule
\textbf{Hyperparameter} & \textbf{Value} \\
\midrule

\multicolumn{2}{@{}l}{\emph{Network architecture}} \\
Shared backbone, fully connected layers & 3 \\
Hidden units per layer & 300 \\
Activation function & ELU \\
Location embedding & Learnable, dimension 100 \\

\addlinespace
\multicolumn{2}{@{}l}{\emph{Horizon and reward}} \\
Decision epoch length $\EpochLength$ & 5 minutes \\
Discount factor $\DiscountRate$ & 1 \\

\addlinespace
\multicolumn{2}{@{}l}{\emph{Optimization and updates}} \\
Number of training days $\Days$ & 500 \\
Minibatch size $|\mathcal{M}_{\text{batch}}|$ & 10 \\
Update frequency $K$ & 1, every decision epoch \\
Replay warm-up threshold & 1{,}000 experiences \\
Optimizer & Adam, learning rate $10^{-3}$ \\
Learning-rate schedule & Adaptive reduce-on-plateau \\
Soft target-update rate $\StepSize$ & $10^{-3}$ \\

\addlinespace
\multicolumn{2}{@{}l}{\emph{Replay and exploration}} \\
Prioritized replay buffer capacity & 5{,}000 experiences \\
Prioritized-replay exponents &
Priority $0.6$, importance sampling annealed $0.4 \rightarrow 1.0$ \\
Matching exploration noise &
Gaussian, std.\ $1$, mean annealed $-10 \rightarrow 0$ \\
Pricing action exploration &
Boltzmann sampling, temperature $1$ \\

\bottomrule
\end{tabular}

\endgroup
\end{table}

\section{Detailed Operational Analysis of the Sensitivity Analysis Experiments}\label{ec:opanalysis}

This section provides the detailed operational analysis underlying the three sensitivity analysis experiments of Section~\ref{results} in the main paper; the sample-deadline experiment is analyzed directly in that section rather than here. For each experiment, the main paper reports the headline service rates (Tables~\ref{table:deadline_table}--\ref{table:locations_table}); here we examine \emph{how} each policy achieves those outcomes by analyzing dispatch frequency, route composition, average travel time, and spatial coverage.

\subsection{Operational Profile under Varying Vehicle Capacity}
Capacity primarily affects how policies balance dispatch frequency, route duration, and volume per trip. Figure~\ref{fig:capacity_tradeoff} presents these three dimensions under the tightest setting ($\Capacity = 10$). Under this constraint, policies that do not adapt their routing structure experience the largest performance losses. This pattern is most evident for \ILPVolume, which falls to 79.07\% at $\Capacity = 10$. As shown in Figure~\ref{fig:capacity_tradeoff}, \ILPVolume maintains low dispatch frequency and long travel times while attempting to maximize volume per trip. With limited capacity, these longer routes cannot collect enough volume to justify their duration, leading to delays and sample expirations. At the opposite extreme, \ILPThroughput dispatches frequently with short travel times but collects little volume per trip, which supports responsiveness but uses vehicle capacity inefficiently. \CombinedNeurADP operates between these extremes: it increases dispatch frequency and shortens routes relative to the consolidated policies, while maintaining substantially higher volume per trip than \ILPThroughput. The ablation variants further clarify the value of the two components. \CentreNeurADP behaves similarly to \ILPVolume, maintaining high consolidation with long travel times and therefore degrading sharply at low capacity. \VehicleNeurADP partially adapts by increasing dispatch frequency and reducing consolidation, but less strongly than \CombinedNeurADP. This suggests that the vehicle value function drives responsiveness and route restructuring, while the centre value function alone tends toward over-consolidation when capacity is tight. As capacity increases to 15 and 20, these differences diminish: policies converge toward more consolidated routing patterns, performance gaps narrow, and \CombinedNeurADP's average volume per trip remains nearly unchanged, indicating that additional capacity is largely underutilized in this tested demand regime.

\begin{figure}[!ht]
\centering
\includegraphics[width=0.99\textwidth]{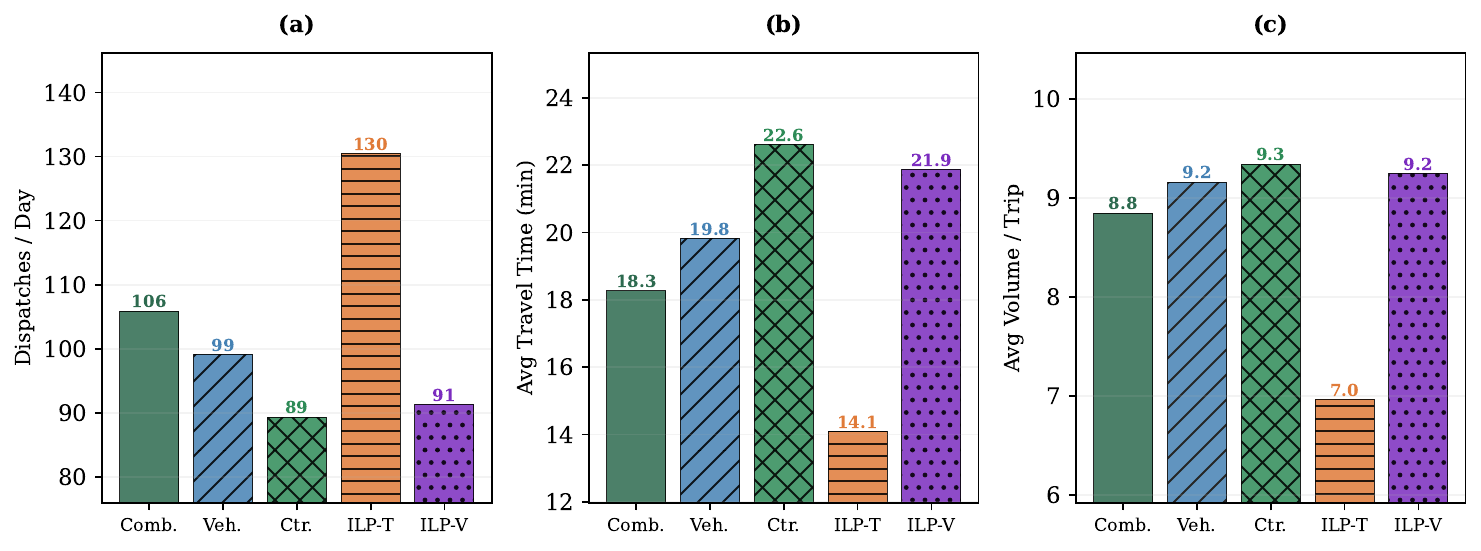}
\caption{Dispatch frequency (a), average travel time (b), and average volume per trip (c) across policies under $\Capacity = 10$.}
\label{fig:capacity_tradeoff}
\end{figure}

\subsection{Operational Profile under Varying Fleet Size}
When vehicles are limited, performance is driven by how effectively each dispatch is utilized. \CombinedNeurADP responds by maximizing the value of each trip, operating with high consolidation and extracting the highest volume per trip among all policies. This reflects a strategy of prioritizing fewer, more productive dispatches when fleet availability is constrained. In contrast, the myopic baselines again fall at opposing extremes. \ILPThroughput fragments the already limited fleet into many short trips, increasing dispatch frequency but collecting insufficient volume per trip. This also leads to severe spatial imbalance, with certain centres receiving disproportionately few visits. As shown in Figure~\ref{fig:agents_balance}, \ILPThroughput exhibits substantially higher variation in centre visits under small fleet sizes, resulting in neglected locations where samples accumulate and expire. At the other extreme, \ILPVolume maintains a highly consolidated strategy with balanced coverage, but its slower circulation prevents it from keeping pace with demand, resulting in the lowest service rate under scarcity. The ablation variants further clarify these dynamics. \CentreNeurADP prioritizes balanced coverage and achieves the lowest variation in centre visits across all fleet sizes, but fails to circulate the limited fleet quickly enough, leading to significant degradation at 3 vehicles. \VehicleNeurADP performs better, adopting a consolidation level similar to \CombinedNeurADP, but lacks the centre-level awareness needed to consistently select the most efficient routes. Together, these results highlight that both components are required: the vehicle value function drives efficient utilization of limited fleet capacity, while the centre value function ensures that this utilization is directed toward the most critical locations. As fleet size increases, these differences diminish. At 5 vehicles, all policies converge toward high performance, with service rates above 95\% and gaps narrowing to approximately 1 percentage point. The additional fleet capacity allows even myopic strategies to compensate for inefficiencies through increased availability, while spatial imbalance is also reduced across all policies (Figure~\ref{fig:agents_balance}). Consistent with earlier findings, these results show that the value of foresighted dispatching is greatest under resource scarcity, where the ability to balance utilization and coverage is most critical.

\begin{figure}[!ht]
\centering
\includegraphics[width=0.609\textwidth]{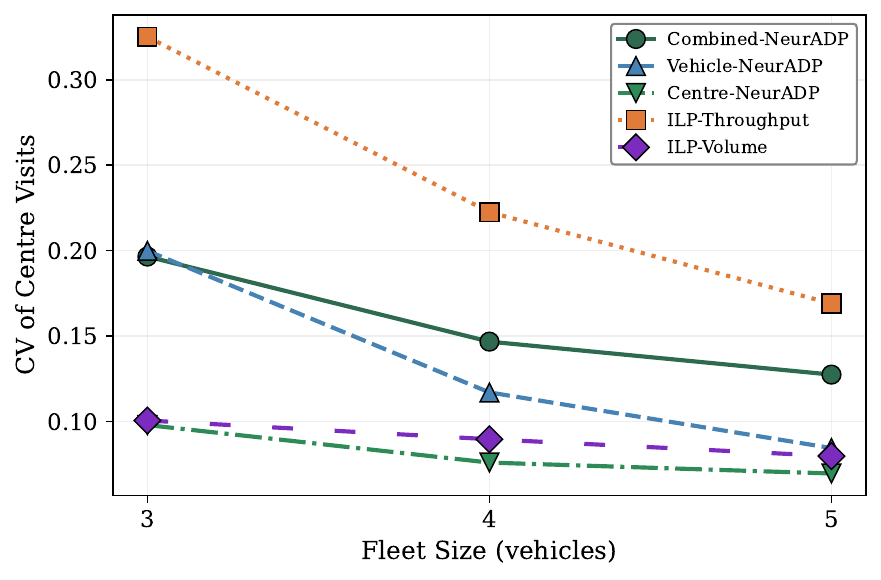}
\caption{Coefficient of variation (CV) of centre visits across policies under varying fleet sizes. Lower values indicate more balanced spatial coverage.}
\label{fig:agents_balance}
\end{figure}

\subsection{Operational Profile under Varying Route Length}
When vehicles are restricted to single-centre trips, routing is trivial and performance depends almost entirely on centre selection. In this setting, \CombinedNeurADP achieves 88.76\%, primarily by maintaining balanced spatial coverage. Notably, \CentreNeurADP slightly outperforms it, reflecting that centre-level information alone is sufficient when route composition is fixed. By contrast, \VehicleNeurADP performs worse due to higher visit imbalance, indicating that without explicit centre-level guidance it fails to prioritize critical locations effectively. Allowing two centres per trip yields a substantial improvement for all policies, as vehicles can begin to consolidate nearby pickups. \CombinedNeurADP benefits most, rising to 92.96\% while maintaining moderate travel times and a high proportion of two-centre routes. In contrast, \ILPThroughput continues to favour single-centre trips even when multi-stop routes are available, limiting its ability to exploit consolidation. The relative ordering of policies aligns with previous experiments, with \CombinedNeurADP outperforming both ablation variants by balancing consolidation with coverage. The most informative setting is $\LocsToVisit = 3$, where the expanded action space introduces both greater flexibility and greater risk. \CombinedNeurADP achieves the highest performance (93.58\%) by selecting a balanced mix of single-, double-, and triple-centre routes, increasing volume per trip while keeping travel times controlled. The ablation variants diverge sharply in this regime. \VehicleNeurADP adopts a similar level of consolidation but exhibits significantly higher imbalance in centre visits, indicating that without centre-level guidance it forms less effective multi-stop routes. In contrast, \CentreNeurADP collapses to 84.74\%, behaving similarly to \ILPVolume by over-emphasizing consolidation. It heavily favors three-centre routes, leading to the longest travel times and reduced responsiveness. \ILPThroughput avoids over-consolidation by continuing to favor short routes, but its fragmentation prevents it from benefiting from the expanded action space. Overall, these results show that as routing complexity increases, the cost of poor trade-offs becomes more pronounced. The vehicle component is required to limit excessive route length, while the centre component ensures adequate spatial coverage. Their combination is therefore most valuable in richer decision spaces, where effective performance requires balancing both considerations.
\end{APPENDICES}

\end{document}